\documentstyle{article}
\begin{document}
\begin{large}

\centerline{\bf TOPOLOGIES ON THE SET OF ALL SUBSPACES}
\centerline{\bf OF A BANACH SPACE AND RELATED QUESTIONS}
\centerline{\bf OF BANACH SPACE GEOMETRY}
\centerline{\bf M.I.OSTROVSKII}

{\bf 1. Introduction.}

For a Banach space $X$ we shall denote  the  set  of  all  closed
subspaces of $X$ by $G(X)$. In some kinds of problems it turned  out  to
be useful to endow $G(X)$ with a topology. The  main  purpose  of  the
present paper is to survey results on two the most common topologies
on $G(X)$.

The organization of this paper is as follows. In section  2  we
introduce some definitions and notation. In  sections  3  and  4  we
introduce two topologies on $G(X)$.  Section  5  is  devoted  to  the
problem  of  comparison  of  these  topologies.  In  section  6   we
investigate the following general problem: How close should  be  the
structure of the subspaces which  are  close  with  respect  to  the
natural metrics, which generate introduced topologies? (It should be
mentioned that  both  introduced  topologies  are  metrizable.)  In
section 7 we survey  those  results  on  introduced  topologies  and
related quantities which were not discussed in  previous  sections.
Here we  also  try  to  describe  known  applications  of introduced
topologies and related quantities. This section is nothing more than
guide to the literature.
If $x$ is a vector of a Banach space $X$ and $A, D$ are subsets of $X$
then we shall denote the value $\inf_{a\in A}||x-a||$  by  dist$(x,A)$  and
the
value $\inf_{a\in A}$dist$(a,D)$ by dist$(A,D)$ . The closed unit ball and the
unit sphere of a Banach  space $X$  are  denoted  by $B(X)$  and $S(X)$
respectively. For a subset $A$ of a  Banach  space $X$  by
$A^{\perp }$, lin($A$),
conv($A$)  and cl($A$)  we  shall   denote,   respectively,   the   set
$\{x^*\in X^*:(\forall x\in A)(x^*(x)=0)\}$, the set of all finite linear
combinations of
vectors of $A$ , the set of all convex combinations of  vectors  of $A$
and the closure of $A$ in the strong topology. For a  subset  $A$  of  a
dual Banach space $X^*$ we shall denote the set
$\{x\in X:(\forall x^*\in A)(x^*(x)=0)\}$
by $A^{\top}$.

Let $Y$ and $Z$ be Banach spaces. For $1\le p\le \infty $ we shall denote
by $Y\oplus _pZ$
the Banach space of  all  pairs $(y,z), y\in Y, z\in Z$,  with  the  norm
$|| (y,z)|| =(|| y|| ^p+|| z|| ^p)^{1/p}
($or $\max \{|| y|| ,|| z|| \}$, if $p=\infty )$. It is  clear  that
all  these  norms  define  the  same  topology.  The   corresponding
topological vector space will be denoted by $Y\oplus Z$.

A closed linear subspace $Y$ of a Banach space $X$ is said to be  a
{\it complemented} subspace of $X$ if there is a bounded  linear  projection
from $X$ onto $Y$ or, what is the same, if there exists a closed  linear
subspace $Z$ of $X$ such that every $x\in X$  can  be  in  a  unique  manner
represented in the form $x=y+z$, where $y\in Y$ and $z\in Z$.
By $G_c(X)$ we  shall
denote the set of all complemented subspaces of $X$.  For $Y\in G_c(X)$  we
shall denote by $\lambda (Y,X)$ the value $\inf \{||P|| :P$ is a projection
of $X$ onto
$Y\}$.  Banach  space  is  said  to  be  {\it injective}  if  its  isomorphic
embeddings into arbitrary Banach space have complemented images.

If $\{X_n\}^{\infty }_{n=1}$ is a sequence of Banach spaces we
define the  direct
sum, of these spaces in the sense of $l_p, 1\le p<\infty $,
 namely $(\sum^{\infty }_{n=1}\oplus X_n)_p$,
as the space of all sequences $x=(x_1,x_2,\ldots
,x_n,\ldots
)$  with $x_n\in X_n$  for
all $n$, for which
$|| x|| =(\sum^{\infty }_{n=1}|| x_n|| ^p)^{1/p}<\infty $.
Similarly, $(\sum^{\infty }_{n=1}\oplus X_n)_{0}$ denotes
the direct sum of $\{X_n\}^{\infty }_{n=1}$ in
the sense of $c_0$, i.e. the space of  all
sequences $x=(x_1,x_{2},\ldots,x_n,\ldots)$ with $x_n\in X_n$
for  all $n$,  for  which
$\lim_n|| x_n|| =0.$ The norm in this direct sum is taken as
$||x|| =\max _n||x_n||$.

A sequence $\{X_n\}^{\infty }_{n=1}$ of closed subspaces of a Banach space
$X$  is
called a {\it Schauder decomposition} of $X$  if  every $x\in X$
has  a  unique
representation of the form $x=\sum^{\infty }_{n=1}x_n$ with
$x_n\in X_n$ for every $n$. In  such
case we write $X=\sum^{\infty }_{n=1}\oplus X_n$. Furthermore,
if $\{Y_n\}^{\infty }_{n=1}$ is  a  sequence  of
subspaces, $Y_n\subset X_n$, then we shall denote
cl(lin$(\cup ^{\infty }_{n=1}Y_n))$  by
$\sum^{\infty }_{n=1}\oplus Y_n$.
It is clear that in this case $\{Y_n\}^{\infty }_{n=1}$ form a
Schauder  decomposition
of $\sum^{\infty }_{n=1}\oplus Y_n$.

Two Banach spaces $Y$ and $Z$ are called {\it isomorphic} if there exists
an invertible operator from $Y$  onto $Z$.  The  {\it Banach-Mazur  distance}
between $Y$ and $Z$ is defined by $d(Y,Z)=\inf ||T|| ||T^{-1}||$, the
infimum  being
taken over all invertible operators from $Y$  onto $Z$ (if $Y$  is  not
isomorphic to $Z$ we put $d(Y,Z)=\infty )$. Let $\{X_n\}$ and $\{Y_n\}$
be two sequences
of Banach spaces. We shall say that they are {\it uniformly isomorphic} if
$\sup _nd(X_n,Y_n)<\infty$.

The identity operator of a Banach space $X$ is denoted by $I_X$ (or
simply by $I$ if $X$ is clear from the context.)
For a mapping $\varphi :A\to B$  by
${\rm im}\varphi $ we denote the set
$\{b\in B:(\exists a\in A)(b=\varphi a)\}$.

Let $T$ be a linear mapping from a linear subspace $L$ of a  Banach
space $X$ into a Banach space $Y$. Then $L$ is called the {\it domain} of
$T$  and
is denoted by $D(T)$. The set $\{x\in D(T): Tx=0\}$ is called the
{\it kernel} of $T$
and is denoted by ker$T$. For Banach spaces $X, Y$ by $L(X,Y)$  we  denote
the space of all bounded linear operators from $X$  into $Y$,  endowed
with the norm $||T|| =\sup \{||Tx|| :x\in B(X)\}$.

The least cardinal $\alpha $ for which there exists a dense subset of $X$
of cardinality $\alpha $ is called the {\it density character} of $X$ and is
denoted by dens$X$.

Referring to formula (3.5) we mean formula (5) from section 3.

{\bf 3. Geometric opening.}

{\bf 3.1. Definition.} Let $Y,Z\in G(X)$.  The  {\it geometric  opening}
(or  simply
{\it opening}, sometimes {\it gap}) between $Y$ and $Z$ is defined to be

$$
\Theta (Y,Z)=\max \{\sup _{y\in S(Y)}\hbox{dist}(y,Z),
\sup _{z\in S(Z)}\hbox{dist}(z,Y)\}
\eqno{(1)}$$
If the sphere of a subspace is empty (it happens when  the  subspace
is $\{0\})$ then the corresponding supremum is set equal to zero.

{\bf 3.2. Remark.} This definition can be used for nonclosed subspaces  as
well. For nonclosed $Y$ and $Z$ we have
$\Theta (Y,Z)=\Theta ($cl$Y$,cl$Z$).

{\bf 3.3.} For $Y,Z\in G(X)$ we let
$$
\Theta _0(Y,Z)=\sup _{y\in S(Y)}\hbox{dist}(y,Z).
$$

{\bf 3.4.} Properties of  the  geometric  opening  are  described  in  the
following theorem.

{\bf Theorem.} {\it Let $X$ be a Banach space and let $Y,Z\in G(X)$. Then}

(a) $0\le \Theta _0(Y,Z)\le 1$ {\it and} $0\le \Theta (Y,Z)\le 1.$

(b) $ \Theta (Y,Z)=\Theta (Z,Y);$

(c) $ \Theta (Y,Z)=0\Rightarrow Y=Z;$

(d) $\Theta _0(Y,Z)=\Theta _0(Z^{\perp },Y^{\perp })$ {\it and therefore}
$$\Theta (Y,Z)=\Theta (Y^{\perp },Z^{\perp }).
\eqno{(2)}$$

(e) {\it If at least one of the subspaces $Y$ and $Z$ is  finite  dimensional
and $\Theta (Y,Z)<1$ then both of them are finite dimensional and}
dim$Y$=dim$Z$.

(f) {\it For every $Y_1,Y_2,Y_3\in G(X)$ the following inequalities hold:}
$$
\Theta _0(Y_1,Y_3)\le \Theta _0(Y_1,Y_2)+
\Theta _0(Y_2,Y_3)+\Theta _0(Y_1,Y_2)\Theta _0(Y_2,Y_3);
$$
$$
\Theta (Y_1,Y_3)\le \Theta (Y_1,Y_2)+\Theta (Y_2,Y_3)+
\Theta (Y_1,Y_2)\Theta (Y_2,Y_3).
\eqno{(3)}$$
  (Inequality (3) is called ``weakened triangle inequality'').

(g) {\it The function $d_g(Y,Z)=\log (1+\Theta (Y,Z))$ is a metric on $G(X)$.
The  set $G(X)$ is a complete metric space with respect to this metric.}

(h) {\it If $X$ is a Hilbert space, then}
$$
\Theta (Y,Z)=|| P_Y-P_Z|| ,
\eqno{(4)}$$
{\it where $P_Y$  and $P_Z$  are  orthogonal  projections  onto $Y$   and $Z$
respectively.}

Proof. Verification of (a), (b) and (c) is immediate.
(d) Recall well-known formulas from the duality theory of Banach
spaces. For $Z\in G(X),\ y\in X$ and $y^*\in X^*$ we have
$$
\hbox{dist}(y,Z)=\sup _{z^*\in S(Z^\perp)}|z^*(y)|;
$$
$$
\hbox{dist}(y^*,Z^{\perp })=\sup _{z\in S(Z)}|y^*(z)|.
$$
Therefore
$$\Theta _0(Y,Z)=\sup _{y\in S(Y)}\hbox {dist}(y,Z)=$$
$$
=\sup _{y\in S(Y),\ z^*\in S(Z^\perp)}|z^*(y)|=
$$
$$
=\sup _{z^*\in S(Z^\perp)}\hbox{dist}(z^*,Y^{\perp })=
\Theta _0(Z^{\perp },Y^{\perp }).
$$

Statement (e) follows immediately from the following lemma.

{\bf 3.5. Lemma.} {\it Let $Y,Z\in G(X),\ Z$ be finite dimensional,}
dim$Y\ge $dim$Z$.  {\it Then
there exists $y\in S(Y)$ such that}
dist$(y,Z)=1.$

Proof. It is clear  that  we  may  suppose  that $Y$  is  finite
dimensional and dim$Y$ =dim$Z$+1. Suppose first that $X$ is strictly convex
i.e. that $|| x_1+x_2|| <||x_1|| +||x_2||$ for every linearly independent
$x_1,x_2\in X$.
It is easy to see that in this case for every $x\in X$  there  exists  a
unique element of best approximation of $x$ by elements of $Z$.  Let  us
denote this element by $A(x)$. The mapping $A$ is, generally speaking,
nonlinear but it is  easy  to  verify  that  it  is  continuous  and
$A(-x)=-A(x)$. We shall use the following topological result.

{\bf 3.6. Theorem.} (K.Borsuk [Bor], see also [DKL]). {\it Let $S^{n-1}$
be the unit
sphere in ${\bf R}^n$ and let $\varphi :S^{n-1}\to {\bf R}^{n-1}$ be
a continuous mapping, for  which
$\varphi (-x)=-\varphi (x)$. Then there exists a point
$x\in S^{n-1}$ such that}
$\varphi (x)=0.$

Since the unit sphere of $Y$ is homeomorphic to $S^{\hbox{dim}Z}$  and
$Z$  is
homeomorphic  to ${\bf R}^{n-1}$  then  this  theorem  is  applicable
in  our
situation and we can find $y\in S(Y)$ such that such that 0 is  the  best
approximation of $y$ by elements of $Z$ i.e. dist$(y,Z)=1.$

Let us turn to the general case. We may assume that $X$=lin$(Y\cup Z)$
and therefore, $X$ is finite dimensional. Let $\{x_{i}\}^n_{i=1}$ be a
basis of $X$
and $\{x^*_i\}$ be its biorthogonal functionals. It is easy to  check  that
for every $k\in {\bf N}$ the norm
$$
||x||_{k}=(||x||^2+(1/k)\sum^n_{i=1}(x^*_i(x))^2)^{1/2}
$$
is strictly convex. Therefore for every $k\in {\bf N}$ we can  find
$y_k\in Y$  such
that $|| y_k|| _k=1$ and dist$_k(y_k,Z)=1.$ Since
$|| y_k|| \le || y_k|| _k=1$ then the sequence
$\{y_k\}^{\infty }_{k=1}$ contains a convergent subsequence.
It is easy to verify that
its  limit  is  the  required  vector.  Lemma  3.5   and   therefore
statement (e) have been proved.

(f) It is clear that we need to prove only the first inequality.
Let $\varepsilon >0$  and $y_1\in S(Y_1)$.  Then  for  some $y_2\in Y_2$
we  have $|| y_1-y_2|| <\Theta _0(Y_1,Y_2)+\varepsilon $,
hence $||y_2||<1+\Theta _0(Y_1,Y_2)+\varepsilon $.
For  some $y_3\in Y_3$  we  have
$||y_2-y_3||<(\Theta _0(Y_2,Y_3)+\varepsilon )|| y_2|| $.
Hence $|| y_1-y_3|| <|| y_1-y_2|| +|| y_2-y_3|| <\Theta _0(Y_1,Y_2)+
\varepsilon +(\Theta _0(Y_2,Y_3)+\varepsilon )
(1+\Theta _0(Y_1,Y_2)+\varepsilon )$. Taking supremum over $y_1\in S(Y_1)$  and
then infimum over $\varepsilon >0$ we obtain the required inequality.

(g) By (a), (b) and (c) the only thing which we need to prove is
the triangle inequality. Let $Y_1,Y_2,Y_3\in G(X)$. By (f) we have
$$
d_g(Y_1,Y_3)=\log (1+\Theta (Y_1,Y_3))\le
\log (1+\Theta (Y_1,Y_2)+\Theta (Y_2,Y_3)+\Theta (Y_1,Y_2)\Theta (Y_2,Y_3))=
$$
$$
\log ((1+\Theta (Y_1,Y_2))(1+\Theta (Y_2,Y_3)))=
\log (1+\Theta (Y_1,Y_2))+\log (1+\Theta (Y_2,Y_3))=
$$
$$
d_g(Y_1,Y_2)+d_g(Y_2,Y_3).
$$

Let $\{Y_n\}^{\infty }_{n=1}\subset G(X)$ be a Cauchy sequence with
respect  to $d_g$.  We
need to show that $\{Y_n\}$ is a convergent sequence. It is sufficient to
prove that $\{Y_n\}$ contains a convergent subsequence. Therefore we  may
assume that $\Theta (Y_n,Y_{n+1})<2^{-n}$. Let us introduce
subset $A\subset S(X)$ as a  set
of limits of all strongly convergent sequences $\{y_n\}^{\infty }_{n=1}$,
for  which
$y_n\in S(Y_n)\ (n\in {\bf N})$. Direct verification shows that $A$ is a
unit sphere of
some closed  subspace $Y_0$  of $X$  and  that
$\lim_{n\to \infty }\Theta (Y_n,Y_0)=0.$ Hence
$\lim_{n\to \infty }d_g(Y_n,Y_0)=0$ and (g) is proved.

(h) We have
$$|| P_Y-P_Z|| =\sup _{x\in S(X)}|| (P_Y-P_Z)x|| =
\sup _{x\in S(X)}|| P_Y(I-P_Z)x-(I-P_Y)P_Zx|| .$$
Hence
$$
|| P_Y-P_Z|| =\sup _{x\in S(X)}(|| P_Y(I-P_Z)x|| ^2+|| (I-P_Y)P_Zx|| ^2)^{1/2}.
\eqno{(5)}$$
Since $||P_Yu||$=dist$(u,Y^{\perp })$
and $|| (I-P_Y)u|| $=dist$(u,Y)$ then we have
$|| P_Y-P_Z|| \le \sup _{x\in S(X)}(\Theta _0(Z^{\perp },Y^{\perp })^2
||(I-P_Z)x||^2+\Theta _0(Z,Y)^2|| P_Zx|| ^2)^{1/2}$.

Therefore, using (2) we obtain
$$
|| P_Y-P_Z|| \le \Theta (Y,Z).
$$
On the other hand, taking  supremum  over $x\in S(Z)$  in  (5),  we
obtain
$$
|| P_Y-P_Z|| \ge \Theta _0(Z,Y),
$$
and taking supremum over $x\in S(Z^{\perp })$ in (5), we obtain
$$
|| P_Y-P_Z|| \ge \Theta _0(Z^{\perp },Y^{\perp }).
$$
These inequalities together with statement (d) imply that
$$
|| P_Y-P_Z|| \ge \Theta (Y,Z).
$$
  Statement (h) is proved.

\centerline{{\bf Notes and remarks}}

{\bf  3.7.} The notion of opening between subspaces of a Hilbert space  was
introduced by M.G.Krein and  M.A.Krasnoselskii [KK].  In [KK]  the
Hilbert space versions of parts (a)--(e) of Theorem 3.4 were proved.
In [KK] opening was introduced by  formula  (1).  The  equality  (4)
appeared in the book due to N.I.Akhiezer and I.M.Glazman [AG].  The
notion of geometric opening between subspaces of a Banach space was
introduced by M.G.Krein, M.A.Krasnoselskii and D.P.Milman in [KKM].
Parts (a)--(e) of Theorem  3.4  were  proved  in  this  paper.  Later
V.M.Tikhomirov [T] rediscovered part (e) of Theorem 3.4.  His  proof
is more complicated. It should be noted that the  authors  of  [KKM]
did not use quantity $\Theta _0(Y,Z)$. This quantity was introduced and
investigated by T.Kato in [K1]. The results of parts (f) and  (g)  of
Theorem 3.4 are due to I.C.Gohberg and A.S.Markus [GM1].

Since the sets $U_{\varepsilon ,Y}=\{Z\in G(X):
\ \Theta (Z,Y)<\varepsilon \}$ form a base of open  sets
of the metric space $(G(X),d_g)$ we shall  say  that  its  topology  is
generated by the geometric opening.
The main reason for introduction of the notion of opening seems
to be  the  following:  this  notion  is  an  appropriate  tool  for
generalization of the  theory  of  Carleman-von  Neumann  of  defect
numbers of Hermitian operators onto the case of arbitrary  operators
in Hilbert and Banach spaces. Now  we  are  going  to  present  this
generalization.

Let $H$ be a complex Hilbert space and $A$ be a linear mapping from
a linear subspace $D(A)\subset H$ into $H$. The operator $A$ is
called  {\it Hermitian}
if $D(A)$ is dense in $H$ and
$$
(\forall x,y\in D(A))((Ax,y)=(x,Ay)).
$$

For  every $\lambda \in {\bf C}$  by $N(\lambda ,A)$  we  shall  denote
the   subspace
$({\rm im}(A-\lambda I))^{\perp }$.

In [KK] and [KKM] the notion of opening was used to  generalize
the following well-known result (see e.g. [Na, \S 14]).

{\bf 3.8. Theorem.} {\it Let A be a Hermitian operator. Then for every  complex
number $\alpha $ from the open upper halfplane we  have
dim$N(\alpha ,A)$=dim$N(i,A)$
and for every complex number $\beta $ from the open lower halfplane we have}
dim$N(\beta ,A)$=dim$N(-i,A)$.

(Here $i=\sqrt {-1}$.) The numbers dim$N(i,A)$ and  dim$N(-i,A)$  are  called
the {\it defect numbers} of $A$.

In   order   to   formulate   the    Krein-Krasnoselskii-Milman
generalization of 3.8 we need the following notion. Let $X$ and $Y$  be
complex Banach spaces. Let $A$ be  a  linear  mapping  from  a  linear
subspace $D(A)\subset X$ into $Y$. A complex number $\alpha $  is
called  a  {\it point  of
regular type} for $A$ if there exists a real number $c(\alpha )>0$ such that
$$
(\forall f\in D(A))(|| (A-\alpha I)f|| \ge c(\alpha )|| f|| ).
$$

For every $\alpha \in {\bf C}$ let us introduce the subspace
$N(\alpha ,A)=({\rm im}(A-\alpha I))^{\perp }\subset Y^*$.

It is easy to see that for every operator $A$ the set  of  points
of regular type for $A$ is an  open  subset  of ${\bf C}$  and  that  complex
numbers with nontrivial imaginary parts are points of  regular  type
for every Hermitian operator. Therefore the following result due  to
M.G.Krein, M.A.Krasnoselskii and D.P.Milman [KKM, p.~111] is a
generalization of 3.8.

{\bf 3.9. Theorem.} {\it Let $G$ be a connected open subset of ${\bf C}$
consisting  of
points of regular type for $A$. Then the  numbers  {\rm dim}$N(\lambda ,A)$
are  the
same for all} $\lambda \in G$.

Proof. It is sufficient to prove that for every point $\alpha \in G$ there
exists a neighbourhood $W\subset G$ such that for every $\lambda \in W$
we have
dim$N(\lambda ,A)$=dim$N(\alpha ,A)$.

Let us show that we may take $W=\{\lambda :|\lambda -\alpha |<
c(\alpha )/4\}\cap G$. By definition  we
have
$$
||(A-\alpha I)f||\ge c(\alpha )|| f|| \hbox{ for every }f\in D(A).
$$
Hence for every $\lambda \in W$ and every $f\in D(A)$ we have
$$
|| (A-\lambda I)f|| \ge (3/4)c(\alpha )|| f|| ;
$$
$$
|| (A-\lambda I)f-(A-\alpha I)f|| =
|\lambda -\alpha ||| f|| <(1/3)|| (A-\lambda I)f|| ;
\eqno{(6)}$$
$$
|| (A-\lambda I)f-(A-\alpha I)f|| <(1/4)|| (A-\alpha I)f|| .
\eqno{(7)}$$
Inequalities (6) and (7) imply that
$$
\Theta ({\rm im}(A-\lambda I),{\rm im}(A-\alpha I))\le 1/3
$$
Using Theorem 3.4 (d) we obtain
$$
\Theta (N(\alpha ,A),N(\lambda ,A))\le 1/3.
$$
By Theorem 3.4 (e) it  follows  that  dim$N(\alpha ,A)$=dim$N(\lambda ,A)$.
The theorem is proved.

{\bf 3.10.}  Although opening $\Theta $ has some properties of  metric
(see  parts
(a), (b) and (c) of Theorem 3.4) for some Banach spaces it is not  a
metric on $G(X)$. More precisely, I.Gohberg and A.S.Markus [GM1]
observed  that $\Theta $ does not satisfy the  triangle  inequality  for  some
Banach spaces.

{\bf Example.} Let $X=l^2_1$ and let $0<a<b\le 1.$
Let us  introduce  subspaces
$Y_1, Y_2, Y_3\in G(X)$ in the following way:
$$
Y_1=\{(x,0):\ x\in {\bf R}\};
$$
$$
Y_2=\{(x,y):\ y=ax,\ x\in {\bf R}\}
$$
$$
Y_3=\{(x,y):\ y=bx,\ x\in {\bf R}\}
$$
It is not hard to see that $\Theta (Y_1,Y_2)=a,\  \Theta (Y_1,Y_3)=b$
and $\Theta (Y_2,Y_3)=
(b-a)/(1+a)$. (The reader advised to draw the picture.) Since
$a+(b-a)/(1+a)=b-a(b-a)/(1+a),$
  and $a(b-a)/(1+a)>0$ then the triangle inequality is not satisfied.

This example shows that the equality in (3) may be attained  for
a triple of pairwise distinct subspaces.

{\bf 3.11.}  In connection with example 3.10 and part (h) of  Theorem  3.4
A.S.Markus proposed the following problem:  Describe  the  class  of
Banach spaces for which the geometric opening satisfies the triangle
inequality.

{\bf 3.12.}  I.Gohberg and A.S.Markus [GM1] suggested to  consider  the
following in some respects more convenient  metric  which  generates
the same topology on $G(X)$ as the geometric opening.

Let $X$ be a  Banach  space, $Y,Z\in G(X)$.  The  {\it spherical  opening}
between $Y$ and $Z$ is defined to be
$$
\Omega (Y,Z)=\max \{\sup _{y\in S(Y)}\hbox{dist}(y,S(Z)),\
\sup _{z\in S(Z)}\hbox{dist}(z,S(Y))\}.
$$
In the case when $Y$ or $Z$ is $\{0\}$, we let $\Omega (Y,Z)=\Theta (Y,Z)$.

The spherical opening satisfies the following  conditions.

a) $\Omega $ is a metric on $G(X)$. (This statement  immediately  follows
from the following observation: $\Omega (Y,Z)$ coincides with the  Hausdorff
distance between $S(Y)$ and $S(Z))$.

b) $\Theta (Y,Z)\le \Omega (Y,Z)\le 2\Theta (Y,Z).$

c) $G(X)$  is  complete  with  respect  to $\Omega$. (This  statement
immediately follows from part (g) of Theorem 3.4.)

{\bf 3.13.}  It is natural to mention one more modification  of  the
geometric opening. This modification was introduced  by  R.Douady  [Do]
and V.I.Gurarii [Gu1].

Let $X$ be a Banach space, $Y,Z\in G(X)$. The {\it ball opening}  between $Y$
and $Z$ is defined to be
$$
\Lambda (Y,Z)=\max \{\sup _{y\in S(Y)}\hbox{dist}(y,B(Z)),\
\sup _{z\in S(Z)}\hbox{dist}(z,B(Y))\}.
$$
In the case when $Y$ or $Z$ is $\{0\}$, we let $\Lambda (Y,Z)=\Theta (Y,Z)$.

It is easy to verify that the analogues of statements (a), (b),
and (c) from 3.12 are valid for the ball opening. Furthermore it
satisfies the inequality
$$
0\le \Lambda (Y,Z)\le 1
$$
and, if $X$ is a Hilbert space then
$$
\Lambda (Y,Z)=||P_Z-P_Y||.
$$
We shall sometimes use the quantities $\Lambda _0(Y,Z)$ and
$\Omega _0(Y,Z)$ which
are defined analogously to $\Theta _0(Y,Z)$.

{\bf 3.14.} It should be mentioned that the analogues  of  the  duality
formula (2)  fail for the spherical and ball openings.

{\bf Example.}   Let $X=l^2_1,\ 0<a<1.$    Let $Y=\{(x,y):\ y=ax,\
x\in {\bf R}\},\
Z=\{(x,0):\ x\in {\bf R}\}$. We have $X^*=l^2_{\infty },\
Y^{\perp }=\{(x,y):\ y=-x/a,\ x\in {\bf R}\},\  Z^{\perp }=\{(0,y):\
y\in {\bf R}\}$.

It is easy to verify that
$$
\Lambda (Y,Z)=\Omega (Y,Z)=2a/(1+a);
$$
$$
\Lambda (Y^{\perp },Z^{\perp })=\Omega (Y^{\perp },Z^{\perp })=a.
$$
(The reader adviced to draw a picture.) But for $0<a<1$ we have
$$
2a/(1+a)\neq a.
$$

{\bf 3.15. Remark.} J.D.Newburgh [New1] introduced another metric on $G(X)$
which generates the same topology on $G(X)$ as the  geometric  opening
does. This metric was investigated and compared with $\Theta $ by  E.Berkson
[Ber].

{\bf3.16.} It is interesting to note that A.L.Brown [Br2] proved that the
assertion of Theorem 3.6 can be easily deduced if  we  suppose  that
the assertion of Lemma 3.5 is true.

{\bf 4. Operator opening.}

{\bf 4.1.} Let $X$ be a Banach space. By $GL(X)$ we  shall  denote  the  group
(with respect to composition) of all invertible linear operators  on
$X$. For $Y,Z\in G(X)$ let
$$
r_0(Y,Z)=\inf \{||C-I||:\ C\in GL(X),\ C(Y)=Z\},
$$
if the set over which  the  infimum  is  taken  is  not  empty,  and
$r_0(Y,Z)=1$ otherwise.

{\bf Definition.} The {\it operator opening} between $Y$ and $Z$ is defined by
$$
r(Y,Z)=\max \{r_0(Y,Z),r_0(Z,Y)\}.
$$

{\bf4.2.} Properties  of  the  operator  opening  are  described  in  the
following theorem.

{\bf Theorem.} {\it Let $X$ be a Banach space and let} $Y,Z\in G(X)$.

(a) $0\le r_0(Y,Z)\le 1$ {\it and hence} $0\le r(Y,Z)\le 1;$

(b) $r(Y,Z)=r(Z,Y);$

(c) $r_0(Y,Z)\ge \Theta _0(Y,Z)$ {\it and hence} $r(Y,Z)\ge \Theta (Y,Z);$

(d) $r_0(Y,Z)<1\Rightarrow r_0(Z,Y)\le r_0(Y,Z)/(1-r_0(Y,Z));$

(e) {\it If $P_Y$ and $P_Z$ are projections with images $Y$ 
and $Z$  respectively,
then} $r(Y,Z)\le ||P_Y-P_Z||$;

(f) $r(Y,Z)=0\Leftrightarrow Y=Z.$

(g) {\it If $X$ is a Hilbert space then for  every $Y,Z\in G(X)$  the  following
equality is valid}
$$
r(Y,Z)=\Theta (Y,Z).
$$

(h) $r_0(Z^{\perp },Y^{\perp })\le r_0(Y,Z).$

(i) {\it For every $Y_1,Y_2,Y_3\in G(X)$ the following inequality is valid:
$$
r_0(Y_1,Y_3)\le r_0(Y_1,Y_2)+r_0(Y_2,Y_3)+r_0(Y_1,Y_2)r_0(Y_2,Y_3)
$$
and hence}
$$
r(Y_1,Y_3)\le r(Y_1,Y_2)+r(Y_2,Y_3)+r(Y_1,Y_2)r(Y_2,Y_3).
$$
(The last inequality is called ``weakened triangle inequality''.)

(j) {\it The function $d_{op}(Y,Z)=\log (1+r(Y,Z))$ is a metric
on $G(X)$. The set
$G(X)$ is complete with respect to metric} $d_{op}$.

Proof. (a) follows from  the  following  observation:  the  set
$\{C\in GL(X):\ C(Y)=Z\}$ with every operator contains all its multiples.

(b) is evident.

(c) If $C_0\in \{C\in GL(X):C(Y)=Z\}$ then
$$
||C_0-I||\ge \sup _{y\in S(Y)}||y-C_0y||\ge \sup _{y\in S(Y)}\hbox{dist}(y,Z).
$$
This inequality implies (c).

(d)  If $C\in GL(X)$  and $C(Y)=Z$  then $C^{-1}\in GL(X)$  and $C^{-1}(Z)$=Y.
Furthermore, we have
$$
||C^{-1}-I||\le ||C-I|| || C^{-1}||\le ||C-I||/\inf _{x\in S(X)}||Cx||\le
||C-I||/(1-||C-I||).
$$
Whence we have (d).

(e) The statement (e) is clear when $||P_Z-P_Y||\ge 1.$ So  we  shall
suppose that $|| P_Z-P_Y|| <1.$ In this case  the  operators $I-(P_Y-P_Z)$  and
$I-(P_Z-P_Y)$ belong to $GL(X)$. Indeed, it is easy  to  verify  that  the
operators $I+\sum^{\infty }_{n=1}(P_Y-P_Z)^n$ and
$I+\sum^{\infty }_{n=1}(P_Z-P_Y)^n$ are  their  inverses.
Therefore $(I-(P_Y-P_Z))X=X$  and  hence $P_Y(I-P_Y+P_Z)X=P_YX$,   therefore
$P_YZ=Y$. Let us denote the operator $I-(P_Z-P_Y)$ by $C$. We  have
$C\in GL(X);\ ||C-I|| =|| P_Z-P_Y||$ and $C(Z)=(I-P_Z+P_Y)Z=Y$.
By  definition  of $r_0(Z,Y)$  we
have $r_0(Z,Y)\le || C-I|| =|| P_Y-P_Z||$.
In  the  same  manner  we  can  estimate
$r_0(Y,Z)$. So (e) is proved.

Statement (f) immediately follows from (c)  and  the  analogous
statement about $\Theta $.

(g) By (e) we have $r(Y,Z)\le || P_Y-P_Z|| $, where $P_Y$  and $P_Z$  are  the
orthogonal projections onto $Y$ and $Z$ respectively. By Theorem 3.4 (h)
$|| P_Y-P_Z|| =\Theta (Y,Z)$.
So we have $r(Y,Z)\le \Theta (Y,Z)$. Comparing this  inequality
with (c) we obtain the desired inequality.

(h)  The  assertion  immediately  follows  from  the  following
observation: if $C\in GL(X)$ and $C(Y)=Z$
then $C^*\in GL(X^*)$ and $C^*(Z^{\perp })=Y^{\perp }$.

(i) The assertion is clear if $r_0(Y_1,Y_2)=1$ or $r_0(Y_2,Y_3)=1.$ So we
shall suppose that it is not the case.

Let $\varepsilon >0$ be arbitrary positive number.
Let $C_1\in GL(X)$ be such that
$C_1(Y_1)=Y_2$ and $|| C_1-I|| <r_0(Y_1,Y_2)+\varepsilon $
and  let $C_2\in GL(X)$  be  such  that
$C_2(Y_2)=Y_3$ and $|| C_2-I|| <r_0(Y_2,Y_3)+\varepsilon $.
Then $C_2C_1\in GL(X),\ C_2C_1(Y_1)=Y_3$ and
$$
|| C_2C_1-I|| =|| (C_2-I)(C_1-I)+(C_1-1)+(C_2-I)|| \le || C_2-I|| || C_1-I|| +
$$
$$
|| C_1-I|| +|| C_2-I|| <
$$
$$
(r_0(Y_2,Y_3)+\varepsilon )
(r_0(Y_1,Y_2)+\varepsilon )+r_0(Y_1,Y_2)+\varepsilon +
r_0(Y_2,Y_3)+\varepsilon .
$$
Since $\varepsilon $ is arbitrary then the required inequality follows.

(j) Since we already proved statements (a), (b) and (f) we need
to verify the triangle inequality and completeness only.

The triangle inequality for $d_{op}$ follows from (i) (see the proof
of Theorem 3.4 (g)).

Let sequence $\{Y_n\}^{\infty }_{n=1}\subset G(X)$ be such that
$$
\lim_{m,n\to \infty }d_{op}(Y_n,Y_{m})=0.
$$
  Since it  is  enough  to  prove  that $\{Y_n\}$  contains  a  convergent
subsequence, we may assume that
$$
r_0(Y_n,Y_{n+1})\le 2^{-n-1}.
$$
Therefore for every $n\in {\bf N}$ there exists an operator
$C_n\in GL(X)$  such
that $C_n(Y_n)=Y_{n+1}$ and $|| C_n-I|| <2^{-n}$.
It is easy to verify that this condition implies the
convergence of the sequence $T_n=\Pi ^n_{k=1}C_k\ (n\in {\bf N})$  in
the uniform topology. Let us denote its limit by $T_0$. It is  easy  to
verify that $T_0\in GL(X)$. Let $Y_0=T_0(Y_1)$.  It  is  easy  to  verify  that
$r_0(Y_k,Y_0)\to 0$ when $k\to \infty $ and, hence,
$\lim_{k\to \infty }d_{op}(Y_k,Y_0)=0.$
We finished the proof of Theorem 4.2.

{\bf4.3. Remark.} It is clear that the sets
$$
V_{\varepsilon ,Y}=\{Z\in G(X):\ r(Y,Z)<\varepsilon \},\
Y\in G(X), \varepsilon >0
$$
form a base of the topology corresponding to the metric $d_{op}$.  So  it
is reasonable to say that this topology is induced by  the  operator
opening.

\centerline{{\bf Notes and remarks}}

{\bf4.4.} Part (g) of Theorem 4.2 implies that in the Hilbert space case
the operator opening is a metric on $G(X)$. However in general case it
is not so: the triangle inequality may fail for it. Indeed, let $X=l^2_1$
and $Y_1,Y_2,Y_3\in G(X)$ be subspaces introduced in example 3.10. It can be
verified that in this case $r(Y_k,Y_{j})=\Theta (Y_k,Y_{j})\ (k,j=1,2,3)$
and so, the
triangle inequality is not satisfied.

{\bf4.5.} The operator  opening  was  introduced  by  J.L.Massera  and
J.J.Sch\"affer [MS1, p.~563]. This concept appeared in a  natural  way
in their investigations of linear differential equations  in  Banach
spaces.  Analogous concept somewhat later was introduced by A.L.Garkavi [Ga].
He used this concept in the theory of best  approximation
in Banach spaces.

{\bf4.6.} All statements of Theorem 4.2 exept statement (h) can be  found
in [Ber], statement (e) was earlier proved by B.Sz-Nagy [Sz1], [Sz2,
p.~132]. Some of them seems to be known to J.L.Massera and  J.J.Sch\"affer
[MS1] (see p.~563). Statement (h) I added since it is a natural
analogue for the corresponding statement about $\Theta $. It seems that  the
example from 4.4 was known to J.L.Massera and  J.J.Sch\"affer  (see
[MS2], \S 13), but it seems that it was not published anywhere.

{\bf 5. Comparison of the topologies induced by  geometric  and  operator}

{\bf openings of subspaces.}

{\bf5.1.} When  geometric  and  operator  openings  were  introduced  the
problem of the comparison of the topologies induced by them arose in
a natural way. By Theorem 4.2(c) the topology induced by the  operator
opening majorize the topology induced by the geometric  opening.
>From Theorem 4.2(g) it follows that in Hilbert space  these  topologies
coincide. In the general case the  following  version  of  this
result is valid.

{\bf 5.2. Theorem.}  {\it Geometric  and  operator  openings  induce  the  same
topology on} $G_{c}(X)$.

This theorem immediately follows from the next result.

{\bf5.3. Proposition.} {\it Let $X$ be a Banach space and $Y\in G_c(X)$.
If $Z\in G(X)$  is
such that $\Omega (Z,Y)<1/\lambda (Y,X)$ then $Z\in G_c(X)$ and}
$$
r_0(Y,Z)\le \lambda (Y,X)\Omega (Z,Y)(1+\lambda (Y,X)-\Omega (Z,Y)
\lambda (Y,X))/(1-\Omega (Z,Y)\lambda (Y,X))
$$

Proof. We shall use  the  following  geometrical  concept.  Let
$Y,Z\in G(X)$. The number
$$
\delta (Z,Y)=\hbox{dist}(S(Z),Y)
$$
is called the {\it inclination} of $Z$ to $Y$.

Let $Y\cap Z=\{0\}$. It is easy to verify that $\delta (Z,Y)>0$ if and only  if
the subspace lin$(Z\cup Y)\subset X$ is closed.
Moreover, if $\delta (Z,Y)>0$ then  there
exists a projection from lin$(Z\cup Y)$ onto $Z$,  whose  kernel  is $Y$  and
whose norm is equal to $1/\delta (Z,Y)$.

By definition of $\lambda (Y,X)$ it follows that  for  every
$\varepsilon >0$  there
exists a projection $P_{Y,\varepsilon }:X\to Y$ such
that $|| P_{Y,\varepsilon }|| <\lambda (Y,X)+\varepsilon $.  Therefore
for $U_{\varepsilon }$=ker$P_{Y,\varepsilon }$ we have
$\delta (Y,U_{\varepsilon })>1/(\lambda (Y,X)+\varepsilon )$.
Using the definition of
$\Omega $ we obtain $\delta (Z,U_{\varepsilon })>\delta (Y,U_{\varepsilon })-
\Omega (Y,Z)>(1-\Omega (Z,Y)(\lambda (Y,X)+\varepsilon ))/
(\lambda (Y,X)+\varepsilon )$.
So, if $\varepsilon $ is small enough then $\delta (Z,U_{\varepsilon })>0.$
Let us show that  if $\varepsilon $  is
small enough then lin$(Z\cup U_{\varepsilon })=X$.
In fact, let us suppose  that  it  is
not the case. Since we may assume without loss of  generality  that
lin$(Z\cup U_{\varepsilon })$  is  closed,  then  we   can
find $x\in S(X)$   such   that
dist$(x$,lin$(Z\cup U_{\varepsilon })>1-\varepsilon $.
Furthermore, we have $x=y+u$,  where $y\in Y, u\in U_{\varepsilon }$
and $|| y|| <\lambda (Y,X)+\varepsilon $.
Therefore   there   exists $z\in Z$   such   that
$|| z-y|| <(\lambda (Y,X)+\varepsilon )(\Omega (Z,Y)+\varepsilon )$.

Hence
$1-\varepsilon <$dist$(x$,lin$(Z\cup U_{\varepsilon }))\le
|| x-(z+u)|| <(\lambda (Y,X)+\varepsilon )(\Omega (Z,Y)+\varepsilon )$.

  It is clear that for $\varepsilon $ small enough this inequality is false.
So we may assume that there exists a projection  of $X$  onto $Z$
whose norm is not greater than $(\lambda (Y,X)+\varepsilon )/
(1-\Omega (Z,Y)(\lambda (Y,X)+\varepsilon ))$  and
whose kernel is $U_{\varepsilon }$.
We shall denote this projection by $P_{Z,\varepsilon }$.  Using
Theorem 4.2(e) we obtain
$$
r(Y,Z)\le \inf _{\varepsilon }|| P_{Z,\varepsilon }-P_{Y,\varepsilon }|| .
$$
  We have
$$
|| P_{Z,\varepsilon }-P_{Y,\varepsilon }||
=\sup _{x\in S(X)}|| P_{Z,\varepsilon }x-P_{Y,\varepsilon }x||
=\sup _{x\in S(X)}|| (P_{Z,\varepsilon }-I)P_{Y,\varepsilon }x|| \le
$$
$$
\sup _{x\in S(X)}\inf _{z\in Z}|| (P_{Z,\varepsilon }-I)
(P_{Y,\varepsilon }x-z)|| \le (|| P_{Z,\varepsilon }|| +1)
|| P_{Y,\varepsilon }|| \Omega (Y,Z).
$$
Taking infimum over $\varepsilon >0$  we  obtain  the  required  inequality.
The proposition is proved.

{\bf5.4.}  V.I.Gurarii  and  A.S.Markus  [GuM]  proved  that  in  general
topologies induced on $G(X)$ by geometric and  operator  openings  are
different. Here are their arguments.
Let $K$ be a complemented subspace of a Banach space $Y$ and  let $L$
be an uncomplemented subspace of a Banach space $Z$ and let $K$ and $L$ be
isomorphic. V.I.Gurarii and A.S.Markus proved that in  such  a  case
the topologies induced on $G(Y\oplus Z)$ by geometric and operator  openings
are different.
\par
In fact, let $T:L\to K$ be an  isomorphism.  Let  us  introduce  the
following family of subspaces of $Y\oplus Z$:
$$
L(\lambda )=\{(\lambda Tx,x):\ x\in L\}\ (\lambda \in {\bf R},\ \lambda >0).
$$
It  is  easy  to  verify  that $\lim_{\lambda \to 0}\Theta (L(\lambda ),L)=0.$
Since $L$   is uncomplemented then in order to prove that
$\lim_{\lambda \to 0}r(L(\lambda ),L)\neq 0$  it  is
sufficient to verify that $L(\lambda )\ (\lambda >0)$ are  complemented.
Let  us  do that. Let $P:Y\to K$ be a projection onto $K$.
Let $P_{\lambda }:Y\oplus Z\to Y\oplus Z$ be defined by
the equality $P_{\lambda }(y,z)=(Py,\lambda ^{-1}T^{-1}Py)$.
It can be directly verified that
$P_{\lambda }$ is a continuous projection onto $L(\lambda )$.

{\bf5.5.} The theory of  complemented  and  uncomplemented  subspaces  is
highly developed now. Comparing results of this theory [BDGJN], [B],
[LT1], [Pel2], [Ros] with the arguments of 5.4  we  can  deduce  the
distinction of the topologies  induced  by  geometric  and  operator
openings for the most part of common Banach spaces  (but  with  such
exeptions as $c_0, l_{\infty }$, Hilbert spaces).

Below  we  shall  describe  another  method  of   proving   the
distinction of the topologies  induced  by  geometric  and  operator
openings. This method works in the cases $X=c_0,l_{\infty }$.
But the  following
problem posed by V.I.Gurarii and A.S.Markus [GuM] remains unsolved.

{\bf Problem.} Does there exist a Banach space which is nonisomorphic
to a Hilbert space but is such that the topologies induced  on $G(X)$
by geometric and operator openings coincide?

At the moment it is known that these topologies  are  different
for Banach spaces which are  not  ``almost  Hilbert''  (in  the  sense
described below).

Recall some definitions. Let us denote by $r_i(t)\ (i\in {\bf N},\ t\in [0,1])$
the Rademacher functions. A Banach space $X$ is said to have {\it type} $p\
(1\le p\le 2)$ if for some  constant $T_p(X)<\infty $  and
for  every  finite  set
$\{x_i\}^n_{i=1}$ of vectors of $X$ the following inequality takes place
$$
(\int^1_0||\sum^n_{i=1}r_i(t)x_i||^2dt)^{1/2}\le
T_p(X)(\sum^n_{i=1}|| x_i|| ^p)^{1/p}.
$$
A Banach space $X$ is said to have {\it cotype} $q\ (2\le q\le \infty )$
if  for  some
constant $C_q(X)<\infty $ and every finite set $\{x_i\}^n_{i=1}$
of vectors  of $X$  the
following inequality takes place:
$$
C_q(X)(\int^1_0|| \sum^n_{i=1}r_i(t)x_i||^2dt)^{1/2}\ge
(\sum^n_{i=1}|| x_i|| ^q)^{1/q}.
$$
Let $p(X)=\sup \{p: X$ has type $p\}$ and $q(X)=\inf \{q: X$ has cotype $q\}$.
It is known that for  infinite  compact $K$  we  have $p(C(K))=1$  and
$q(C(K))=\infty $. If a measure $\mu $ is distinct from the
atomic  measure  with
finite  set  of  atoms  and $1\le r<\infty $
then $p(L_{r}(\mu ))=\min \{r,2\}$   and
$q(L_{r}(\mu ))=\max \{r,2\}$. S.Kwapien [Kw] proved that
if a Banach space  has
type 2 and cotype 2 then it is isomorphic to a Hilbert  space.  This
result allows us to consider a Banach space $X$ with $q(X)=p(X)=2$ as an
``almost Hilbert'' one.

Extensive information on type  and  cotype  (including  results
mentioned above) may be found in [LT2], [MiS] and [Pis1].

{\bf5.6 Theorem.} {\it If a Banach space $X$ is such that
$p(X)\neq 2$ or $q(X)\neq 2,$ then
the topology induced on $G(X)$ by the  operator  opening  is  strictly
stronger than the topology induced by the geometric opening.}

We divide the proof of the theorem into two  parts.  The  first,
Proposition 5.7, is a ``pasting together''  of  infinite-dimensional
subspaces $\{W_j\}^{\infty }_{j=0}$  for  which  the   following  
conditions   are
satisfied:
$$\lim_{j\to \infty }r(W_0,W_j)\neq 0,
\eqno{(1)}$$
$$\lim_{j\to \infty }\Theta (W_0,W_j)=0.
\eqno{(2)}$$
from finite-dimensional ``pieces''. The second, Proposition 5.8, is
a construction of suitable  finite-dimensional  ``pieces''  in
spaces satisfying the condition of the theorem.

{\bf5.7. Proposition.} {\it Let a Banach space $X$ be such that  for
some $\gamma >1,$
every $\varepsilon >0$ and every $Y\in G(X)$ of finite codimension in $X$
there  exist
finite-dimensional subspaces $Z_1$ and $Z_2$ in $G(Y)$ such that
$\Theta (Z_1,Z_2)<\varepsilon $
and $d(Z_1,Z_2)>\gamma $. Then the topology induced by the operator opening on
$G(X)$ is strictly stronger than the topology induced by the geometric
opening on} $G(X)$.

Proof. We construct a sequence of subspaces of $X$ that satisfy
conditions (1) and (2). Let the numbers
$\{\delta _i\}^{\infty }_{i=1}$ be such that $\delta _i>0$
and $\Delta =\Pi ^{\infty }_{i=1}(1+2\delta _i)<\infty $.
We  take $Y_1=X,\ \varepsilon _1=1.$  According  to   the
condition of the proposition, finite dimensional $Z^1_1, Z^1_2\subset Y_1$
with
$\Theta (Z^1_1,Z^1_2)<\varepsilon _1,\  d(Z^1_1,Z^1_2)>\gamma $ can  be  found.
We  consider  the $\delta _1$-net
$\{y_i\}^{n(1)}_{i=1}$ on $S($lin$(Z^1_1\cup Z^1_2))$
and functionals $\{y^*_i\}^{n(1)}_{i=1}\subset S(X^*)$ such that
$y^*_i(y_i)=1.$ We let $Y_2=(\{y^*_i\}^{n(1)}_{i=1})^{\top}$
(intersection of the  kernels  of
functionals $y^*_i$), $\varepsilon _2=1/2.$   For   them
in   turn   we    find
finite-dimensional $Z^2_1,Z^2_2\subset Y_2$, such that
$\Theta (Z^2_1,Z^2_2)<\varepsilon _2,\ d(Z^2_1,Z^2_2)>\gamma $. We
consider  the $\delta _2$-net $\{y_i\}^{n(2)}_{i=1}$
of $S($lin$(\cup ^2_{i=1}(Z^i_1\cup Z^i_2)))$   and
functionals $\{y^*_i\}^{n(2)}_{i=1}$  such  that $y^*_i(y_i)=1.$
Further,   we   take
$Y_3=(\{y^*_i\}^{n(2)}_{i=1})^{\top},\ \varepsilon _3=1/3,$  etc.
We  shall  investigate  the  space
$V=$cl(lin$(\cup ^{\infty }_{i=1}Z^i))$, where $Z^i$=lin$(Z^i_1\cup Z^i_2)$.
Repeating the arguments from [LT1, p.~4] we can prove  that  spaces
$Z^i$ form a Schauder decomposition of $V$.
Let $\pi _{k,n},\ k\le n$,  be
operators on $V$ defined  by  the  equations
$\pi _{k,n}(\sum^{\infty }_{i=1}z^i)=\sum^n_{i=k}z^i$.  By
arguments analogous to those in [LT1, p.~5], we can  show  that  the
operators $\pi _{k,n}$  are  bounded   and
$\sup _{k,n}|| \pi _{k,n}|| \le 2\Delta <\infty $.   Therefore
${\bf [}z{\bf ]}=\sup _{k,n}|| \pi _{k,n}z|| $ is an equivalent norm on $V$.
We extend it onto  the
whole $X$ using the following construction due to A.Pelczynski  [Pel1]
(Proposition 1). Let $B_1$ be the unit ball of $V$ in the  new  norm  and
let $\alpha >0$ be such that $(\alpha B(X)\cap V)\subset B_1$, where $B(X)$
is the unit ball of $X$
in the original norm. Let us introduce new norm on $X$ as a  Minkowski
functional of cl(conv$(\alpha B(X)\cup B_1))$.
It is easy to verify that this norm
coincides with {\bf [$\cdot$]} on $V$.
So we may denote it also by {\bf [$\cdot$]}. It is easy
to verify  that  this  norm  coincides  with  the  original  on
$Z^i,\ i=1,2,\ldots$, and in addition has the following property:
$$
(\forall m(2)>m(1)>n(1)>n(2))(\forall (z^i)^{m(2)}_{i=n(2)}, z^i\in Z^i)
$$
$$
({\bf [}\sum^{m(2)}_{i=n(2)}z^i{\bf ]\ge [}\sum^{m(1)}_{i=n(1)}z^i{\bf ]}).
\eqno{(3)}$$
It is clear that it is sufficient to prove  the  distinction  of
the topologies for $X$ with the norm {\bf [$\cdot$]}.
The  desired  sequence  of
subspaces will be the following:
$$
W_0=\sum^{\infty }_{i=1}\oplus Z^i_1;\
W_j=(\sum^{}_{i\neq j}\oplus Z^i_1)\oplus Z^j_2.
$$
In order to prove (2) it is sufficient to show that
$\Theta (W_j,W_0)<\varepsilon _j$.
Let $z\in S(W_j),\ z=z^j_2+\sum_{i\neq j}z^i_1$.
By (3) we  have ${\bf [}z^j_2{\bf ]}\le 1.$  Using
the fact that a new norm coincides on $Z^j$ with  the  original,  and
consequently the inequality $\Theta (Z^j_1,Z^j_2)<\varepsilon _j$
is preserved,  we  find  a
vector $z^j_1\in Z^j_1$ such that ${\bf [}z^j_1-z^j_2{\bf ]}<\varepsilon _j$.
For the vector $\tilde{z}=\sum^{\infty }_{i=1}z^i_1$ we have
$\tilde{z}\in W_0, {\bf [}\tilde{z}-z{\bf ]}<\varepsilon _j$.
Analogous arguments can be carried  out  also  for
$\tilde{z}\in S(W_0)$. Thus (2) is proved.

We shall prove (1) by contradiction. Let $\{\phi _j\}^{\infty }_{j=1}$
be  a  sequence
of bounded linear operators on $X$ that satisfy the conditions:
$$
\hbox{I)}\  {\bf [}I_X-\phi _j{\bf ]\to }0 ;
$$
$$
\hbox{II)}\  \phi _j(W_j)=W_0.
$$
We introduce operator $\tau_j:Z^j_2\to Z^j_1$ as the  restriction
of  operator
$\pi _{j,j}\phi _j$ to the space $Z^j_2$. For $z^j_2\in Z^j_2$ we have
$\phi _j(z^j_2)=\tau_jz^j_2+\sum^{}_{i\neq j}w^i_1$  for
some $w^i_1\in Z^i_1$. Consequently, by (2) we have
$$
{\bf [}z^j_2-\tau_jz^j_2{\bf ]}\le{\bf [}(z^j_2-\tau_jz^j_2)+
\sum^{}_{i\neq j}w^i_1{\bf ]}={\bf [}z^j_2-\phi _jz^j_2
{\bf ]}\le{\bf [}z^j_2{\bf ][}I_X-\phi _j{\bf ]}.
$$
We obtain ${\bf [}z^j_2-\tau_jz^j_2{\bf ]}\to 0$ uniformly over
$z^j_2\in S(Z^j_2)$. Therefore, starting
from some $j$, the operators $\tau_j$ are isomorphisms of $Z^j_2$
into $Z^j_1$, and,
in addition ${\bf [}\tau_j{\bf ][}\tau^{-1}_j{\bf ]}\to 1$  when
$j\to \infty $.  This  contradicts  inequality
$d(Z^j_1,Z^j_2)>\gamma >1,$ which is preserved also in the new  norm,
since  it
coincides on $Z^j$ with the original. The proposition is proved.

{\bf5.8. Proposition.} {\it If a Banach space $X$ is  such  that $q(X)\neq 2$
or $p(X)\neq 2,$ then for every $\gamma >1$  and  every
$\varepsilon >0$  there  exist  finite
dimensional  subspaces $Z_1,Z_2\in G(X)$  such   that $d(Z_1,Z_2)>\gamma $
and} $\Theta (Z_1,Z_2)\le \varepsilon $.

In the proof of this proposition  we  shall  use  the  following
construction. Let $X$ be a Banach  space  and  let $Y\in G(X)$.  Let  us
denote by $\vartheta $ be the quotient mapping $X\to X/Y$. In the  space
$X\oplus _1(X/Y)$
we introduce the following subspaces $G_0=Y\oplus _1(X/Y)$ and
$G_{\varepsilon }=\{(\varepsilon x,\vartheta x):\
x\in X\},\ 1>\varepsilon >0.$

{\bf5.9. Lemma.} {\it The equalities
$d(G_{\varepsilon },X)\le (1+\varepsilon )/\varepsilon $ and
$\Theta (G_0,G_{\varepsilon })\le \varepsilon $  take
place.}

Proof. The first inequality follows from the fact  that  the  operator
$\tau:X\to G_{\varepsilon }$, defined by the equation
$\tau x=(\varepsilon x,\vartheta x)$, satisfies  inequalities
$|| \tau|| \le 1+\varepsilon , || \tau^{-1}|| \le 1/\varepsilon $.
We shall prove the second inequality. Let
$u\in G_{\varepsilon },\ u=(\varepsilon x,\vartheta x)$.
Then $|| u|| =\varepsilon || x|| +|| \vartheta x|| $.
By the fact that $\vartheta $  is  a  quotient mapping, for  every
$\delta >0$  we  can  find $y_{\delta }\in $Ker$\vartheta =Y$  such  that
$|| x-y_{\delta }|| <|| \vartheta x|| +\delta $. We introduce a vector
$v_{\delta }(u)=(\varepsilon y_{\delta },\ \vartheta x)\in G_0$.
For $v_{\delta }(u)$ we
have $|| u-v_{\delta }(u)|| <(|| \vartheta x|| +\delta )\varepsilon $.
Therefore,
$$
\sup _{u\in S(G_{\varepsilon})}\hbox{dist}(u,G_0)\le
\sup \{\inf _{\delta >0}(|| u-v_{\delta }(u)|| /|| u|| ):\
u\in G_{\varepsilon },\ u\neq 0\}\le
$$
$$
\sup \{\inf _{\delta >0}\varepsilon (|| \vartheta x|| +
\delta )/(\varepsilon || x|| +|| \vartheta x|| ):\
u=(\varepsilon x,\vartheta x)\neq 0\}\le \varepsilon .
$$
Let $v\in G_0,\ v=(y,z)$. We take $x_{\delta }\in X$ such that
$z=\vartheta x_{\delta }$ and $|| x_{\delta }|| \le || z|| +\delta $.
We introduce vectors $u_{\delta }(v)=(y+\varepsilon
x_{\delta },\ z)\in G_{\varepsilon }$. We have
$$
\sup _{v\in S(G_0)}\hbox{dist}(v,\ G_{\varepsilon })\le \sup \{
\inf _{\delta >0}(|| v-u_{\delta }(v)|| /|| v|| ):\ v\in G_0, v\neq 0\}\le
$$
$$
\sup \{\inf _{\delta >0}\varepsilon (|| z|| +\delta )/(|| y|| +|| z|| ):\
0\neq v=(y,\ z)\in G_0\}\le \varepsilon .
$$
Lemma 5.9 is proved.

{\bf5.10.} Proof of Proposition 5.8. Let $1<p=p(X)<2,\  q=p/(p-1)$. In  [Ros,
p.~286] it was proved that there exists a sequence of  spaces $\{X_n\}$,
dim$X_n=n$, that satisfy the following conditions:

1) $X_n$ is a subspace of $l^{m}_q$ for some $m=m(n)\in {\bf N}$.

2) The subspaces $X^{\perp }_n\subset l^{m}_p$ are such that for any
sequence $\{Y_n\}$  of
subspaces of $l_p$ which are uniformly isomorphic to $\{X^{\perp }_n\}$
we have
$$
\lim_{n\to \infty }\lambda (Y_n,l_p)=\infty .
$$

3) Spaces $\{X^*_n\}$ are uniformly isomorphic to some subspaces $\{W_n\}$
of $l_p$.

We introduce the spaces
$$
L_n=l^{m(n)}_p\oplus _1(l^{m(n)}_p/(X_n^{\perp })),
\eqno{(4)}$$
  where sequences $\{X_n\}^{\infty }_{n=1}$ and
$\{m(n)\}^{\infty }_{n=1}$ are such that conditions 1,  2
and 3 are satisfied. We note that the second term  in  (4)  is
isometric to $X^*_n$, and therefore from condition 3 and  the  Maurey-Pisier
theorem [MiS, p.~85] it follows that a sequence of  subspaces $M_n\subset X,
\ n=1,2,\ldots$ can be found such that $\sup _nd(M_n,L_n)<\infty $.
Therefore  it  is
sufficient to construct the desired $Z_1$ and $Z_2$ in some of $L_n$.  After
using Lemma 5.9 for $X=l^{m(n)}_p,\ Y=X_n^{\perp }$, we find $G^n_0,\
G^n_{\varepsilon }\subset L_n$  such  that
$\Theta (G^n_0,G^n_{\varepsilon })\le \varepsilon ;\
 d(l^{m(n)}_p,G^n_{\varepsilon })\le (1+\varepsilon )/\varepsilon $.
>From condition 2 and  the  equation
$G^n_0=X_n^{\perp }\oplus _1(l^{m(n)}_p/(X_n^{\perp }))$
we conclude that $d(l^{m(n)}_p,G^n_0)\to \infty $  when $n\to \infty $.
Since $d(G^n_0,G^n_{\varepsilon })\ge d(l^{m(n)}_p,G^n_0)/d(l^{m(n)}_p,\
G^n_{\varepsilon })$, then by  choosing $n$
sufficiently large, we get $d(G^n_0,G^n_{\varepsilon })>\gamma $.
Therefore for $Z_1$  and $Z_2$  we  can
take, respectively, $G^n_{\varepsilon }$ and $G^n_0$.
\par
In the case $p(X)=1$ we use the fact that  the  above-constructed
$Z_1$ and $Z_2$ can be embedded in $X$ with the help of M.I.Kadets'  theorem
[Kad1] (see also [MiS, p.~50]).

We consider now the case $q=q(X)>2.$ We shall use the  fact  that
for $p=q/(q-1)\ (p=1$ in the case $q=\infty )$ we have the following  assertion
[MiS, pp.~21, 23]: $\alpha >0$ can be found such that for each
$n\in {\bf N}$ in $l^n_p, a$
subspace $Y_n$  can  be  found  with  dim$Y_n=[\alpha n]$ (integral  part)  and
$d(Y_n.l^{[\alpha n]}_2)<2.$ In the same time, if subspaces
$Z_n\subset l^n_q$  are  such  that
dim$Z_n/n^{2/q}\to \infty $ when $n\to \infty$  (dim$Z_n/\ln (n)\to \infty $
in the case $q=\infty )$, then
$$
d(Z_n,l^{t(n)}_2)\to \infty ,\hbox{ when }n\to \infty ,
\eqno{(5)}$$
  where $t(n)$=dim$Z_n$.  From  the  first  statement  it  follows  that
subspaces $X_n$ of codimension $[\alpha n]$ can be  found  in
$l^n_q$  such  that
$d(l^n_q/X_n,l^{[\alpha n]}_2)<2.$

We  consider  spaces $L_n=l^n_q\oplus _1(l^n_q/X_n)$.
Using  the  Maurey-Pisier
theorem [MiS, p.~85] and Dvoretzky's theorem [MiS, p.~24], we  get
that  subspaces $M_n$ can be found in $X$  such  that
$\sup _nd(L_n,M_n)<\infty $.
Therefore it is sufficient to construct the desired $Z_1$ and $Z_2$  in
some of $L_n'$s. For this we use the lemma with $X=l^n_q,\ Y=X_n$. We obtain
spaces $G^n_{\varepsilon }, G^n_0\subset L_n$ for  which
$d(G^n_{\varepsilon },l^n_q)\le (1+\varepsilon )/\varepsilon,\
\Theta (G^n_{\varepsilon },G^n_0)\le \varepsilon $.  Since
$G^n_0=X_n\oplus _1(l^n_q/X_n)$ has a subspace
$l^n_q/X_n$  for  which $d(l^n_q/X_n,l^{[\alpha n]}_2)<2,$
then from (5) it  follows  that $d(G^n_0,l^n_q)\to \infty $
when $n\to \infty $.  We  argue
further in the same way as in the first case. The  proposition  is
proved.

{\bf 5.11.} The theorem follows from Propositions 5.7 and 5.8  since,
as it is easy to see, if $Y\subset X$ and codim$Y<\infty $, then
$p(X)=p(Y)$ and $q(X)=q(Y)$.

\centerline{{\bf Notes and remarks.}}

{\bf5.12.} Theorem 5.2 is due to E.Berkson [Ber]. Theorem  5.6  and
Propositions 5.7 and 5.8 are due to the author [O7]. After  publication
of [O7] the author learned that Lemma 5.9 was  already  known  to
A.Douady [D] (see pp. 15--16). In [O4] the ``localized'' version of 5.4
was developed. Using this version it can be proved that in order  to
solve in negative the Gurarii-Markus problem from 5.5 it  is  sufficient
to prove the following inversion of  B.Maurey's  theorem  from
[Ma]: If a Banach space $X$ does not have type 2 then there  exists  a
sequence  of  finite  dimensional   subspaces $X(n)\subset X$
such   that
$\sup _nd(X(n),l^{\dim (X(n))}_2)<\infty $ and
$\sup _n\lambda (X(n),X)=\infty $.

\centerline {\bf 6. Community of properties of subspaces which are close}
\centerline {\bf with respect to the opening.}

{\bf6.1.} Subspaces which are close with respect to the operator  opening
are close with respect to the Banach-Mazur distance.

{\bf Proposition}. {\it If $r_0(Y,Z)<1,$ then}
$$
d(Y,Z)\le (1+r_0(Y,Z))/(1-r_0(Y,Z)).
$$

Indeed, let $T\in GL(X)$ be such  that $T(Y)=Z$  and $||T-I||\le
r_0(Y,Z)+\varepsilon <1.$
Then $T^{-1}=I+\sum^\infty_{i=1}(I-T)^{i}$. Therefore
$$
||T^{-1}||\le 1+\sum^{\infty }_{i=1}(r_0(Y,Z)+
\varepsilon )^{i}=1/(1-r_0(Y,Z)-\varepsilon ).
$$
Whence it follows the desired inequality.

{\bf6.2.} The following proposition may be considered as an  inversion of
statement from 6.1.

{\bf Proposition.}{\it Let $Y$ and $Z$ be isomorphic
Banach  spaces.  Then  for
every $\varepsilon >0$ there
exists a Banach space $X$ and isometric embeddings  of
$Y$ and $Z$ into $X$ such that for their images (which we still denote  by
$Y$ and $Z$) the following inequalities are satisfied:}
$$
r_0(Y,Z)\le d(Y,Z)+\varepsilon -1;
\eqno{(1)}$$
$$
\Omega (Y,Z)\le d(Y,Z)+\varepsilon -1
\eqno{(2)}$$

Proof.  Let $U:Y\to Z$  be  an  isomorphism  such   that $||U||=1$   and
$||U^{-1}||\le d(Y,Z)+\varepsilon $. On the direct
sum $Y\oplus Z$  we  introduce  the  following
seminorm:
$$
p(y,z)=\max\{||z+Uy||,\sup\{|z^{*}(z)+U^*z^*(y)/||U^*z^*|| |:z^*\in S(Z^*)\}\}.
$$

This seminorm generates a norm on the quotient  of $Y\oplus Z$  by  the
zero space of this seminorm. We shall denote the completion  of  the
corresponding normed space by $X$. The spaces $Y$  and $Z$  isometrically
embed into $X$ in a natural way.

Let us introduce an operator $T$ on $Y\oplus Z$ by the equality:
$$
T(y,z)=(-U^{-1}z, 2z+Uy).
$$
We have
$$
p((T-I)(y,z))=p(-U^{-1}z-y,z+Uy)=
$$
$$
\max\{||(z+Uy)+U(-U^{-1}z-y)||,\sup\{|z^*(z+Uy)+U^*z^*(-U^{-1}z-y)/||U^*z^*|| |:
$$
$$
z^*\in S(Z^*)\}\}\le||z+Uy||\sup\{|1-1/||U^*z^*|| |:z^*\in S(Z^*)\}\le
$$
$$
p(y,z)(||(U^*)^{-1}||-1)\le p(y,z)(d(Y,Z)+\varepsilon -1).
$$
This inequality implies that $T$ induces continuous operator on $X$
for which $T(Y)=Z$ and $||T-I||\le d(Y,Z)+\varepsilon-1.$
Thus, we have inequality (1).

In order to prove inequality (2) it is sufficient to verify that
for every $y\in S(Y)$
we have $p(y,-Uy/||Uy||)\le d(Y,Z)+\varepsilon -1.$ We have
$$
p(y,-Uy/||Uy||)=
$$
$$
\max\{||Uy||(1/||Uy||-1);\sup\{|-z^*(Uy)/||Uy||+U^*z^*(y)/||U^*z^*|| | :\
z^*\in S(Z^*)\}\}\le
$$
$$
\max\{(||U^{-1}||-1);\sup\{|z^*(Uy)\left|\right|1/||U^*z^*||-1/||Uy|| |:\
z^*\in S(Z^*)\}\}\le
$$
$$d(Y,Z)+\varepsilon -1.$$

The proposition is proved.

{\bf6.3.} By Propositions 6.1 and 6.2 the class of properties  which  are
common for subspaces which are close with respect  to  the  operator
opening coincides with the class of properties which are common  for
spaces which are close with respect to  the  Banach-Mazur  distance.
The situation with the analogous problem for the  geometric  opening
turns out to be quite different. No smallness of the $\Theta (Y,Z)$  implies
that $Y$ and $Z$ are isomorphic. This assertion may be  deduced  from
Lemma 5.9. Indeed, it is well-known that for $X=l_1$  there  exists  a
subspace $Y\in G(X)$ such that $X/Y$ is isometric to $l_{2}$ [LT1, p.~108].
Let
$G_{\varepsilon }\  (\varepsilon >0)$ and
$G_0$ be the subspaces introduced  before  Lemma  5.9.  By
Lemma 5.9 all $G_{\varepsilon }\ (\varepsilon >0)$ are
isomorphic to $l_1$ and $\Theta (G_{\varepsilon },G_0)\to 0$
when $\varepsilon \to 0.$
But the space $G_0=Y\oplus l_{2}$ is not isomorphic to $l_1$ [LT1, p.~54].
This arguments prove that the set of all subspaces in $G(X)$ which
are isomorphic to the given one ($Y\oplus l_{2}$ in the example)
may  be  even
non-open in $G(X)$ in the topology induced by the geometric opening.

{\bf6.4.} In this connection it seems natural to introduce the  following
definitions.

By {\it property} we shall mean a subclass in the class of all  Banach
spaces. When $X$ satisfies property $P$ we shall write $X\in P$.

A property $P$ will be called {\it open} if for every Banach space $X$ the
subset $P\cap G(X)\subset G(X)$ is open in the topology induced
by the opening $\Theta $.

A property $P$ is called {\it stable} if
there exists a number $\alpha >0$  such
that for every Banach space $X$
and  every $Y,Z\in G(X)$,  if $Y\in P$  and
$\Theta (Y,Z)<\alpha $, then $Z\in P$. The
least upper bound of  numbers $\alpha $  for  which
this statement is  true  is  denoted  by $s(P)$  and  is  called  the
{\it stability exponent} of $P$.

A property $P$ will be called {\it extendedly stable} if there exists  a
number $\alpha >0$ such that for
every Banach space $X$ and every $Y,Z\in G(X)$, if
$Y\in P$ and $\Theta_0(Z,Y)<\alpha $, then $Z\in $P.
The least upper bound of numbers $\alpha $  for
which this statement is true is denoted by $es(P)$ and is  called  the
{\it extended stability exponent} of $P$.

{\bf 6.5. Remark.} {\it It  is  clear  that  stable  properties  are  open  and
extendedly stable properties are stable.}

{\bf 6.6. Proposition.} {\it If property $P$ is open,
then for  every $Y\in P$  there
exists a number $\alpha >0$ such that for every
isometric  embedding $U:Y\to X$
and every $Z\in G(X)$, if $\Theta (Z,UY)<\alpha $, then $Z\in P$.}

Proof. Let us suppose the contrary. Then there  exists  a  space
$Y\in P$, spaces $X(i)\ (i\in {\bf N})$, isometric embeddings
$U_i:Y\to X(i)$ and subspaces
$Z(i)\in G(X(i))$ such that
$\Theta (U_{i}Y,Z(i))\to 0$ when $i\to \infty $.

Let us introduce a  space $X=(\sum^{\infty }_{i=1}\oplus X(i))_1$
and  a  subspace $K\subset X$
consisting of those sequences $\{x(i)\}^{\infty }_{i=1}$
for  which $x(i)\in {\rm im}U_i$  and
$\sum^{\infty }_{i=1}U^{-1}_{i}x(i)=0.$ Then the restrictions
of the quotient mapping $\varphi :X\to X/K$
to the subspaces $X(i)$  are  isometries.
For  every $i,j\in {\bf N}$  we  have
$\varphi U_{i}Y=\varphi U_{j}$Y. Let us denote
this subspace of $X/K$ by $V$.  We  have $V\in P,\
\varphi (Z(i))\not\in P$ and
$\lim_{i\to \infty }\Theta (\varphi (Z(i)),V)=0.$
This  contradicts  to  the  fact
that $P$ is open. The proposition is proved.

{\bf6.7. Proposition.} {\it Stable properties are isomorphic invariants.}

Proof. Let $P$ be a stable property  and $s(P)$  be  its  stability
exponent.  Using  formula  (2)  with $\varepsilon =s(P)/2$
we  obtain  that  if
$d(Y,Z)<s(P)/2+1$ and $Y\in P$ then $Z\in $P.

Let $Y\in P,\ Z$ be isomorphic to $Y$ and $T:Y\to Z$ be an  isomorphism.  Let
us introduce the following family of equivalent norms on $Y$:
$$
||y||_{t}=(1-t)||y||_Y+t||Ty||_Z,\ t\in [0,1].
$$
 It is clear that $Y(0)$ is isometric to $Y,\ Y(1)$ is isometric to $Z$, and
the function $t\to d(Y,Y(t))$ is continuous. Hence we can find a  set  of
points $t_1,t_{2},\ldots,t_n\in [0,1]$ such that
$$
d(Y,Y(t_1))<s(P)/2+1;\ d(Y(t_n),Z)<s(P)/2+1;
$$
$$
d(Y(t_k),Y(t_{k+1}))<s(P)/2+1,\ (k=1,2,\ldots ,n-1).
$$
 By this and the observation made at the begining  of  the  proof  it
follows that $Z\in P$. The proposition is proved.

{\bf6.8.} An open property need not be an isomorphic invariant.

{\bf Example.} $P=\{X:(\exists x_1,x_2\in B(X))(||x_1+x_{2}||>3/2)\
\&(||x_1-x_{2}||>3/2)\}$.  It  is
not hard to show that this property is open. On the other  hand,  by
the parallellogram identity it
follows that $l_2\not\in P$. In the  same time,
the space $l^{2}_1\oplus _1l_{2}$ is
isomorphic to $l_{2}$ and $l^{2}_1\oplus _1l_{2}\in P$.

{\bf6.9.} Set-theoretic operations (union, intersection, complement)  may
be introduced for properties in a natural way. It  is  not  hard  to
verify the following assertions:

{\bf Proposition.} a) {\it The intersection of finite  collection  of  open
properties is open. The union of every class of open  properties  is
open.}

b) {\it Let $\{P_{\alpha }:\alpha \in A\}$  be
some  class  of  stable  properties   and
$\inf _{\alpha \in A}s(P_{\alpha })>0,$
then $\cap _{\alpha \in A}P_{\alpha }$
and $\cup _{\alpha \in A}P_{\alpha }$  are  stable   properties.
Furthermore,
$s(\cap _{\alpha \in A}P_{\alpha })>\inf _{\alpha \in A}s(P_{\alpha })$
and $s(\cup _{\alpha \in A}P_{\alpha })>\inf _{\alpha \in A}s(P_{\alpha })$.
Analogous assertion is valid for extendedly stable properties.}

c) {\it The complement of the stable property is stable.}

As a rule we shall not formulate  consequences  of  our  results
which can be obtained by immediate application of this proposition.

{\bf6.10. Example.} By Theorem 3.4 (e) it follows that for every  natural
number $n$ the class of all $n$-dimensional subspaces is stable and  its
stability exponent is 1. By Lemma 3.5 the class of the spaces  whose
dimension is not greater than $n$ is extendedly stable.  It  is  clear
that its complement is not extendedly stable.

{\bf6.11.} The following result is an immediate consequence of Lemma  5.9
and Proposition 6.7.

{\bf Proposition.} a) {\it Let property $P$ be  such  that  for  some  Banach
spaces $X$ and $W$ the following conditions are satisfied:}

1) $X\in P$.

2) $X$ {\it can be decomposed into a direct sum} $X=Y\oplus Z$.

3) $W$ {\it contains a subspace $W_0$, which is isomorphic to $Y$  and  such
that the quotient $W/W_0$ is isomorphic to} $Z$.

4) $W$ {\it is not isomorphic to any space from} $P$.

{\it Then $P$ is not open.}

b) {\it If there exists $X\in P$ such  that  for  some $Y\in G(X)$  the  space
$Y\oplus _1(X/Y)$ is not in $P$, then $P$ is not stable.}

c) {\it If $P$ is such that for some Banach spaces $X$ and $W$  conditions
{\rm 1--3 {\it of} (a)} are satisfied and $W$ contains a subspace which is
not  in $P$, then $P$ is not extendedly stable.}

d) {\it If there exists $X\in P$ such  that  for  some $Y\in G(X)$  the  space
$Y\oplus _1(X/Y)$ contains a subspace which is
not  in $P$,  then $P$  is  not
extendedly stable.}

{\bf6.12.} Part (a) of proposition 6.11 shows that if property $P$  is  an
isomorphic invariant and is such that for $Y,Z\in P$
we have $Y\oplus Z\in P$  then
the negative solution of the so-called ``three space  problem''  is  a
sufficient condition for $P$ to be non-open.

Recall that the ``three space problem'' for property $P$  is  the
following problem: Let $X$ be a Banach space
such that for some $Y\in G(X)$
we have $Y\in P$ and $X/Y\in P$. Does it follow that $X\in P$?

This problem was investigated for different properties  by  many
authors. In the present context results of  the  negative  character
are of interest. Such results can be found  in [CG1], [CG2], [C],
[ELP], [JR], [KP], [Lu], [O8], [OP].

It is natural to ask the following question: Let property $P$  be
non-open. Is it true that the ``three space problem'' for $P$  has  the
negative answer? In order  to  avoid  trivial  situations  we  shall
additionally suppose that $P$ is closed with respect to  formation  of
direct sums and is isomorphic invariant.

It turned out that even under this restriction the answer to the
posed question is negative. The corresponding example  is  given  by
the introduced in 5.5 class of ``almost Hilbert'' spaces. In [ELP]  it
was proved that  the  ``three  space  problem''  for  this  class  has
positive solution. On the other hand we shall see  later  that  this
class is not open.

It should be noted that by the known results on type and  cotype
[Pis2] it follows that the class of ``almost Hilbert'' spaces does not
satisfy the condition of part (b) of Proposition 6.11.

{\bf6.13. Example.} Proposition 5.3 immediately implies that the class of
injective Banach spaces is open. Let us
consider $X=l_{\infty }$ and let $Y\in G(X)$
be isometric to $l_{2}$. By part (b) of Proposition 6.11 it follows  that
the class of injective spaces is not stable.

{\bf6.14.} Now we turn to another method of establishing  unstability  of
classes of Banach spaces. This method is quite  general:  every
unstable class satisfies its conditions. But sometimes it is not clear
how to apply this method to concrete classes.

Let $Y$  and $Z$  be  Banach  spaces  and  let $T:S(Y)\to S(Z)$   and
$D:S(Y^*)\to S(Z^*)$ be  surjective  mappings.  Let  us  introduce  on  the
algebraic sum $Y\oplus Z$ the following seminorm:

$$
p(y,z)=\sup\{|y^*(y)-(Dy^*)(z)|:y^*\in S(Y^*)\}.
$$
Seminorm $p$ generates norm on the
quotient of $Y\oplus Z$ by  the  zero-space
of $p$. We denote the completion of this normed space by $X$.

By properties of $D$ it follows that $Y$ and $Z$ are  isometric  with
their natural images in $X$. (This images we shall still denote  by $Y$
and $Z$.) We have
$$
\Omega _0(Y,Z)=\sup \{{\rm dist}(y,S(Z)):y\in S(Y)\}\le
\sup \{||y-Ty||:y\in S(Y)\}=
$$
$$
\sup \{|y^*(y)-(Dy^*)(Ty)|:\ y\in S(Y),\ y^*\in S(Y^*)\}.
$$
$$
\Omega _0(Z,Y)=\sup \{{\rm dist}(z,S(Y)):z\in S(Z)\}\le
\sup \{||z-y||:\ y\in T^{-1}z,\ z\in S(Z)\}=
$$
$$
\sup \{||y-Ty||:\ y\in S(Y)\}.
$$
Hence $\Omega (Y,Z)$ is not greater than the following quantity.
$$
\sup \{|y^*(y)-(Dy^*)(Ty)|:\ y\in S(Y),\ y^*\in S(Y^*)\}.
\eqno{(3)}$$
Let  us  introduce  the  quantity $k(Y,Z)$  as  the  infimum  of
quantities  (3)  over  all  surjective  mappings $T:S(Y)\to S(Z)$   and
$D:S(Y^*)\to S(Z^*)$. We have
$$
\inf _{X,U,V}\Omega (UY,VZ)\le k(Y,Z),
\eqno{(4)}$$
 where the infimum is taken  over  all  Banach  spaces $X$  containing
isometric copies of $Y$ and $Z$ and over all isometric embeddings $U:Y\to X$
and $V:Z\to X$. It turns out [O9] that for certain $a>0$ (e.g., we may  let
$a=1/20$) we have
$$
\inf _{X,U,V}\Omega (UY,VZ)\ge ak(Y,Z).
\eqno{(5)}$$
Using inequalities (4) and (5) and Proposition 6.6 we obtain:

{\bf Proposition.} (a) {\it Property $P$ is non-open if and only if for some
$Y\in P$ and some Banach spaces $Z_n\not\in P\ (n\in {\bf N})$ we have}
$$
\lim_{n\to \infty }k(Y,Z_n)=0.
$$
(b) {\it Property $P$ is unstable if and only if for  some $Y_n\in P$  and  some
Banach spaces $Z_n\not\in P\ (n\in {\bf N})$ we have}
$$
\lim_{n\to \infty }k(Y_n,Z_n)=0.
$$

{\bf6.15.} M.I.Kadets [Kad2] applied construction from 6.14 to $Y=l_{2}$  and
$Z=l_{p}$, where $1<p<2.$ Considering mappings $T:l_{2}\to l_{p}$
and $D:l_{2}\to l_{q}$
(where $q$  is  such  that $1/q+1/p=1$) defined by the equations
$$
T(\{x_{i}\}^{\infty }_{i=1})=
\{|x_{i}|^{2/p}{\rm sign}(x_{i})\}^{\infty }_{i=1};
$$
$$
D(\{x_{i}\}^{\infty }_{i=1})=\{|x_i|^{2/q}{\rm sign}
(x_{i})\}^{\infty }_{i=1},
$$
he obtained the following estimate:
$$
k(l_{2},l_{p})\le 2(2/p-1).
\eqno{(6)}$$

By proposition 6.14 (a) this inequality implies that  the  class  of
spaces isomorphic to $l_{2}$ is non-open.

Using analogous estimates the present author [O1]  proved  that
for every $1<p<\infty $ the class of spaces isomorphic to $l_{p}$ is non-open.

{\bf6.16. Corollary.} {\it The class of ``almost Hilbert'' spaces
is non-open.}

Indeed, since for $r<2$ we have $p(l_{r})=r<2$ (see  5.5)  it  follows
that for every $r<2$ the space $l_{r}$ is not ``almost  Hilbert''.  Comparing
this assertion with Proposition 6.14 and estimate (6) we obtain  the
required result.

{\bf6.17.} Let us turn to stable properties. We are going to describe a
method of finding extendedly stable properties.

Let $\Gamma $ be a set, $l_1(\Gamma )$ be the corresponding
Banach  space.  (I.e.
the space of functions $f:\Gamma \to {\bf R}$  with  countable  support,
denoted  by
supp$f$ and such that $\sum_{\gamma \in {\rm supp}f}
|f(\gamma )|<\infty $. The norm on $l_1(\Gamma )$  is  defined
as $||f||=\sum_{\gamma \in {\rm supp}f}|f(\gamma )|.$) Let $X$ be  a  Banach
space.  By $l_{\infty }(\Gamma ,X)$  we
denote the space of functions $x:\Gamma \to X$ such
that $\sup _{\gamma \in \Gamma }||x(\gamma )||_{X}<\infty $,  with
the norm $||x||=\sup _{\gamma \in \Gamma }||x(\gamma )||_{X}$.

Let $A$ be a subset of the unit sphere of $l_1(\Gamma )$.
For every $a\in A$
we introduce a linear operator from $l_{\infty }(\Gamma ,X)$  into $X$
defined  in  the following way:
$$
x\to \sum_{\gamma \in \Gamma }a(\gamma )x(\gamma ).
$$
This operator will be also denoted by $a$. It is clear that  the  norm
of this operator equals 1.

{\bf6.18. Definition.} By the {\it index of $A$ in} $X$ we mean
the supremum $h(X,A)$
of  those $\delta $  for  which  there  exists
$x\in S(l_{\infty }(\Gamma ,X))$   such   that
$\inf _{a\in A}||a(x)||\ge \delta $.

Many common properties of Banach  spaces  can  be  described  in
terms of introduced indices.

{\bf6.19. Example.} Recall the characterization  of  reflexivity  due  to
R.C.James [J1] and D.P.Milman-V.D.Milman [MM].

Let $X$  be  a  Banach  space.  The  following   assertions   are
equivalent:

a) $X$ is nonreflexive.

b) For some $\varepsilon >0$ there exists a sequence
$\{x_{i}\}^{\infty }_{i=1}\subset B(X)$ such that
$$
(\forall n\in {\bf N})({\rm dist}({\rm conv}\{x_1,x_{2},\ldots
,x_n\},{\rm conv}\{x_{n+1},\ldots \})\ge \varepsilon )
\eqno{(7)}$$

c) For every $1>\varepsilon >0$ there exists a sequence
$\{x_{i}\}^{\infty }_{i=1}\subset B(X)$ such  that
(7) is satisfied.

In order to restate these results in terms of introduced indices
let us define $R\subset S(l_1)$ as a set of all vectors of the form
$$
(a_1,\ldots ,a_n,-a_{n+1},\ldots ,-a_{m},0,\ldots ),
$$
where $n<m$ are arbitrary natural numbers and numbers $a_{i}\ (i=1,\ldots ,m)$
are such that $a_{i}\ge 0,\ \sum^n_{i=1}a_{i}=1/2;\
\sum^{m}_{i=n+1}a_{i}=1/2.$

It is easy to see that the formulated above characterization  of
reflexivity can be restated as following:

A Banach space $X$ is reflexive if and only if $h(X,R)=0.$ If $X$  is
nonreflexive, then $h(X,R)\ge 1/2.$

{\bf6.20. Example.} Let $\alpha $ be  some  uncountable
cardinal.  M.~G.~Krein,
M.~A.~Krasnoselskii and D.~P.~Milman [KKM, p.~98] proved that the following
conditions are equivalent:

(a) dens$X\ge \alpha $.

(b) For some $\varepsilon >0$ the unit ball $B(X)$ contains a
subset of cardinality
$\alpha $ such that the distance between each two elements of
this subset is not less than $\varepsilon $.

(c) For every $0<\varepsilon <1$ the unit ball $B(X)$ contains a
subset  of  cardinality $\alpha $ such that the distance between each
two  elements  of  this subset is not less than $\varepsilon $.

This result also may be reformulated  in  terms  of  introduced
indices. Let set $\Gamma $ be such that card$\Gamma =\alpha $.
Let us denote by $D(\alpha )$  the subset of $S(l_1(\Gamma ))$
consisting of all vectors with two-point  support
for which one of the values is $1/2$ and the other is  $(-1/2)$.  It  is
easy to see that  the  mentioned  characterization  of  spaces  with
density character not  less  than $\alpha $  can  be  reformulated  in  the
following way.

The density character of a Banach space $X$ is less than $\alpha $ if and
only if $h(X,D(\alpha ))=0.$ If dens$X\ge \alpha $,
then $h(x,D(\alpha ))\ge 1/2.$

{\bf6.21. Proposition.} {\it Let $X$ be a Banach space, $Y,Z\in G(X)$.
Then for every
set $\Gamma $ and every $A\subset S(l_1(\Gamma ))$ the following
inequalities take place:}
$$
\Lambda _0(Y,Z)\ge h(Y,A)-h(Z,A)
\eqno{(8)}$$
$$
\Theta_0(Y,Z)(1+h(Z,A))\ge h(Y,A)-h(Z,A)
\eqno{(9)}$$

Proof of (8). Let $\varepsilon >0$ be arbitrary. Let vector
$y\in S(l_{\infty }(\Gamma ,Y))$  be
such that $\inf _{a\in A}||a(y)||\ge h(Y,A)-\varepsilon $.
Using  definition  of $\Lambda _0$  we  find
$z\in B(l_{\infty }(\Gamma ,Z))$ such that
$$
(\forall \gamma \in \Gamma )(||z(\gamma )-y(\gamma )||<
\Lambda _0(Y,Z)+\varepsilon ).
$$
 Therefore\qquad $\inf _{a\in A}||a(z)||\ge h(Y,A)-
\Lambda _0(Y,Z)-2\varepsilon $.            Hence
$h(Z,A)\ge h(Y,A)-\Lambda _0(Y,Z)$.

Proof of (9). Let $\varepsilon >0$ be arbitrary.
Let vector $y\in S(l_{\infty }(\Gamma ,Y))$  be
such that $\inf _{a\in A}||a(y)||\ge h(Y,A)-\varepsilon $.
Using  definition  of $\Theta_0$  we  find
$z\in l_{\infty }(\Gamma ,Z)$ such that
$$
(\forall \gamma \in \Gamma )(||z(\gamma )-y(\gamma )||<\Theta_0(Y,Z)+
\varepsilon ).
$$
 Therefore
$$
\inf _{a\in A}||a(z)||\ge h(Y,A)-\Theta_0(Y,Z)-2\varepsilon .
$$
 On the other hand
$$
||z||\le \sup _{\gamma }||y(\gamma )||+\Theta_0(Y,Z)+\varepsilon .
$$
 Hence
$$
h(Z,A)\ge (h(Y,A)-\Theta_0(Y,Z))/(1+\Theta_0(Y,Z))
$$
 It is easy to see that this  inequality is equivalent  to  (9).  The
proposition is proved.

Applying this proposition to examples 6.19 and  6.20  we  obtain
the following consequences.

{\bf 6.22. Corollary.} {\it Let $P$ be a class of reflexive  spaces. 
Then $P$  is
extendedly stable and} $es(P)\ge 1/2.$

{\bf6.23. Corollary.} {\it Let $\alpha $ be uncountable ordinal and
let $P$ be  a  class
of Banach spaces whose density character is less than $\alpha $. Then $P$  is
extendedly stable and} $es(P)\ge 1/2.$

{\bf 6.24.} It turns out that many  known  isomorphic  invariants  can  be
described in the same manner as  properties  in  examples  6.19  and
6.20. In this connection it is natural to  introduce  the  following
definition.

{\bf Definition.} Class $P$ of Banach spaces is said to be $l_1${\it -property},
if there exist a set $\Gamma $ and a subset $A\subset S(l_1(\Gamma ))$
such  that  for  some
positive $\delta >0$ the following assertions are equivalent:

 (a) $X\not\in P.$

(b) $h(X,A)>0.$

(c) $h(X,A)\ge \delta.$

 By Proposition 6.21  every $l_1$-property  is  extendedly  stable  and
extended stability exponent is not less than $\delta $.

{\bf6.25. Definition.} Let $P$ be an $l_1$-property.
The supremum of  those $\delta $
for which there exist a set $\Gamma $ and a
subset $A\subset S(l_1(\Gamma ))$ such that  the
conditions (a), (b) and (c) of  definition  6.24  are  satisfied  is
called the $l_1$-{\it exponent} of $P$ and is denoted by $e_1(P)$.

It is clear that $es(P)\ge e_1(P)$.

{\bf 6.26.} It seems that $e_1(P)$ can be less than $es(P) ($see  e.g.  example
6.27.4 below). But now I haven't proof of this assertion.

{\bf 6.27.} Here is the list of known $l_1$-properties and the estimates  for
their $l_1$-exponents and extended stability exponents.

{\bf 6.27.1.} Reflexivity. By 6.19  it  follows  that  reflexivity  is  an
$l_1$-property and that its $l_1$-exponent  is  not  less  than  1/2.  The
precise  values  of  the $l_1$-exponent  and  the  extended  stability
exponent of reflexivity seems to be unknown.

{\bf 6.27.2.} Let $\alpha $ be an uncountable cardinal. The class  of 
all  Banach
spaces, whose density character is less than $\alpha $
is an $l_1$-property  by
the result of M.~G.~Krein,
M.~A.~Krasnoselskii and D.~P.~Milman  mentioned
in 6.20, and its $l_1$-exponent is not less than 1/2. The $l_1$-exponents
and the extended stability exponents of those classes for which $\alpha $ is
not greater than the cardinality of continuum are equal to 1/2  (see
6.43). For those $\alpha $  which  are  greater  than  the  cardinality  of
continuum  the  values  of $l_1$-exponents  and   extended   stability
exponents seem to be unknown.

{\bf 6.27.3.} The class of all finite dimensional  spaces.  One  can  show
that this class is an $l_1$-property by  consideration  of  the  subset
$A\subset S(l_1)$ consisting of all vectors with two-point support  for  which
one of the values is 1/2 and the other is (-1/2). Using this set  it
is not hard to verify that the $l_1$-exponent  of  the  class  of  all
finite dimensional spaces is not less than 1/2.  Its  precise  value
seems to be unknown.  On  the  other  hand  by  Theorem  3.4(e)  the
extended stability exponent of this class equals 1.

{\bf 6.27.4.} Let $n\in {\bf N}$. The class of all finite dimensional
Banach  spaces,
whose dimension is not greater than $n$ is an $l_1$-property. It  can  be
shown by consideration of $A=S(l^{n+1}_1)$. Using notion of Auerbach system
(see [LT1, p.~16]) it can be shown that the $l_1$-exponent of the class
of spaces whose dimension is not greater than $n$, is  not  less  than
$1/(n+1)$. It is  known  that  this  estimate  is  not  precise.  This
assertion can be derived, for  example,  from  the  results  on  the
estimates  of  Banach-Mazur  distances  between $l^n_1$  and   arbitrary
$n$-dimensional space (see [Sz]  and  [SzT]).  Even  the  order  (with
respect to $n)$ of these $l_1$-exponents seems to be unknown. By  Theorem
3.4(e) the extended stability exponents of these classes  are  equal
to 1.

{\bf 6.27.5.} $B$-convexity.

{\bf Definition.} A Banach space $X$ is said to be $B${\it -convex} if
$$
\lim_{n\to \infty }
\inf \{d(X_n,l^n_1): X_n\hbox {is an}\ n-
\hbox {dimensional subspace of}\ X\}=\infty .
$$

It was shown by D.P.Giesy [Gi] (Lemmas I.4 and I.6),  see  also
[LT2, p.~62], that for every non-$B$-convex Banach space $X$,
every $n\in {\bf N}$
and every $\varepsilon >0$ there exists a subspace $X_n\subset X$ such
that $d(X_n,l^n_1)<1+\varepsilon $.

Let $A\subset S(l_1)$ be the set of all sequences,  which  for
some $n\in {\bf N}$ have the following form:
$$
(0,\ldots ,0,a_1,\ldots ,a_n,0,\ldots ),
$$
where $a_1$ is preceded by $n(n-1)/2$ zeros.

It is clear that the mentioned  result  due  to  Giesy  can  be
reformulated in the following way.

Let $X$ be a Banach space. The  following  three  assertions  are
equivalent.

(a) $X$ is non-$B$-convex.

(b) $h(X,A)>0.$

(c) $h(X,A)=1.$

Hence $B$-convexity is an $l_1$-property and its $l_1$-exponent and its
extended stability exponent are equal to 1.

{\bf6.27.6.} Super-reflexivity.

{\bf Definition.}  A  Banach  space $X$  is  said   to   be   {\it finitely
representable} in a Banach space $Y$ if  for  each  finite  dimensional
subspace $X_n$ of $X$ and each positive number $\varepsilon $ there exists
a  subspace
$Y_n$ of $Y$ such that $d(X_n,Y_n)<1+\varepsilon$. A Banach  space
$X$  is  said  to  be
{\it super-reflexive}  if   every   Banach   space   which   is   finitely
representable in $X$ is reflexive.

The  class  of  super-reflexive  spaces  has  many   equivalent
descriptions. For our purposes the most interesting is the following
one (see [Be2, p.~236, 237, 270]).

For a Banach space $X$ the following properties are equivalent.

 (a) $X$ is not super-reflexive.

 (b) For some $\varepsilon >0$ and every $n\in {\bf N}$ there
is a sequence $x^n_1,\ldots ,x^n_n$ in  the
unit ball of $X$ such that for every $k\in {\bf N}, 1\le k\le n$,
$$
{\rm dist}({\rm conv}(x^n_1,\ldots
,x^n_k),{\rm conv}(x^n_{k+1},\ldots
,x^n_n))>\varepsilon .
\eqno{(10)}$$

 (c) For every $0<\varepsilon <2$ and every $n\in {\bf N}$ there is a
sequence $x^n_1,\ldots ,x^n_n$  in
the unit ball of $X$ such that inequality (10) is satisfied for  every
$k\in {\bf N}, 1\le k\le n$.

Let $A\subset S(l_1)$ be the set of all sequences,  which  for
some $n\in {\bf N}$ have the following form:
$$(0,\ldots ,0,a_1,\ldots ,a_n,0,\ldots ),$$
where for some $k\in {\bf N},\ 1\le k\le n$, we
have $a_{i}\le 0$  for $i\le k$  and
$a_{i}\ge 0$ for $i\ge k+1,$; $\sum^k_{i=1}a_{i}=-1/2$
and $a_1$ is preceded by $n(n-1)/2$ zeros.

It is clear  that  the  mentioned  characterization  of  super-reflexivity
can be reformulated in the following way.

Let $X$ be a Banach space. The  following  three  assertions  are
equivalent.

(a) $X$ is not super-reflexive.

(b) $h(X,A)>0.$

(c) $h(X,A)=1.$

Hence super-reflexivity is an $l_1$-property and  its $l_1$-exponent
and its extended stability exponent are equal to 1.

{\bf6.27.7.} Class of spaces which do not contain  isomorphic  copies  of
$l_1$. R.C.James [J2] (see also [LT1, p.~97]) proved that if  a  Banach
space $X$ contains a subspace isomorphic to $l_1$, then,
for  every $\varepsilon >0$
there exists a subspace $Y\subset X$ for which
$$
d(Y,l_1)<1+\varepsilon .
$$
Let $A$ be the set of all finitely non-zero vectors  from $S(l_1)$.
It is clear that mentioned result due to James can  be  reformulated
in the following way.

Let $X$ be a Banach space. The  following  three  assertions  are
equivalent.

 (a) $X$ contains an isomorphic copy of $l_1$.

(b) $ h(X,A)>0.$

(c) $ h(X,A)=1.$

Hence the class of spaces which do not contain  isomorphic  copies
of $l_1$ is an $l_1$-property and
its $l_1$-exponent  and  its  extended
stability exponent are equal to 1.

{\bf6.27.8.} Alternate-signs Banach-Saks property.

{\bf Definition.} A Banach space $X$ is said  to  have  {\it alternate-signs
Banach-Saks property}  ($ABS$)  if  every
bounded  sequence $\{x_n\}^\infty _{n=1}\subset X$
contains a subsequence $\{x_{n(i)}\}^{\infty }_{i=1}$  such  that
the  alternate-signs Cesaro  means $n^{-1}\sum^n_{k=1}(-1)^kx_{n(k)}$  are
convergent  in   the   strong topology.

B.Beauzamy [Be1, p.~362] (see also [BL]) obtained the following
characterization of $ABS$. Let $X$ be a Banach
space. Then the following
assertions are equivalent.

 (a) $X\not\in ABS.$

 (b) There exist $\varepsilon >0$ and a sequence $\{x_n\}\subset B(X)$
such that for all $k\in {\bf N}$, if $k\le n(1)\le \ldots
\le n(k)\ (n(i)\in {\bf N})$, then for all scalars $c_1,\ldots
,c_k$,
$$
||\sum^k_{i=1}c_ix_{n(i)}||\ge \varepsilon
(\sum^k_{i=1}|c_{i}|)
\eqno{(11)}$$

(c) For every $0<\varepsilon <1,$ there exists a sequence
$\{x_n\}\subset B(X)$ such that for all $k\in {\bf N}$,
if $k\le n(1)\le \ldots \le n(k)\ (n(i)\in {\bf N})$,  then  for
all   scalars $c_1,\ldots
,c_k$, inequality (11) is satisfied.

Let $A\subset S(l_1)$ be the set of all finitely non-zero  vectors,  for
which the least element of support is not less than its cardinality.
It is clear that  the  mentioned  result  due  to  Beauzamy  can  be
reformulated in the following way.

Let $X$ be a Banach space. The  following  three  assertions  are
equivalent.

 (a) $X\not\in ABS.$

(b) $h(X,A)>0.$

(c) $h(X,A)=1.$

Hence $ABS$ is an $l_1$-property
and its $l_1$-exponent and its extended stability
exponent are equal to 1.

{\bf6.27.9.} Banach-Saks property.

{\bf Definition.} A Banach space $X$ is said to have  {\it Banach-Saks
property} ($BS$) if every  bounded  sequence
$\{x_n\}^{\infty }_{n=1}\subset X$  contains  a  subsequence
$\{x_{n(i)}\}^{\infty }_{i=1}$ such that the Cesaro means
$n^{-1}\sum^n_{k=1}x_{n(k)}$ are convergent in the strong topology.

Let $A\subset S(l_1)$ be the set  consisting  of  all  finitely  non-zero
vectors satisfying the following  conditions.  If
supp$a=\{n(1),\ldots ,n(k)\}$ and $n(1)\le \ldots \le n(k)$, then:

(a) $k\le n(1).$

 (b) There  exists $j\le k$  such  that $a_{n(1)},\ldots ,a_{n(j)}\ge0;\
a_{n(j+1)},\ldots ,a_{n(k)}\le0$ and
$\sum^{j}_{i=1}a_{n(i)}=-\sum^k_{i=j+1}a_{n(i)}=1/2.$

Using the descriptions of reflexivity (6.19), $ABS$  (6.27.8)  it
is not hard to prove that the following statements are equivalent.

 (a) $X\not\in BS.$

(b) $h(X,A)>0.$

(c) $h(X,A)\ge 1/2.$

Hence $BS$ is an $l_1$-property and its $l_1$-exponent and its extended
stability exponent are not less than 1/2. Their precise values  seem
to be unknown.

{\bf6.28.} The list of $l_1$-properties presented above is by no means
complete. Later (6.38) we shall describe methods  of  constructing  new
$l_1$-properties from the known ones.

It should be mentioned that every two $l_1$-properties  from  6.27
are different. Some of the  corresponding  verifications  are  quite
nontrivial (see books [Be2], [BL], [Du] and references in 6.40).

{\bf6.29.} Now we  shall  describe  two  ways  of  obtention  new  stable
properties from the known ones.

Let $P$ be some property of  Banach  spaces.  The  class  of  all
Banach spaces for  which $X^{**}/X\in P$  is  denoted  by $P^{co}$.
J.~Alvarez,
T.~Alvarez and M.~Gonzalez [AAG] proved the following general result.

{\bf Theorem.} {\it If property $P$
is extendedly stable then $P^{co}$  is  also
extendedly stable and $es(P^{co})\ge es(P)/2.$ The  analogous  statement  is
valid  for  stable  properties.  If  property $P$  is  an  isomorphic
invariant and is open then $P^{co}$ is also open.}

In order to prove this theorem we need the following result due
to M.~Valdivia. (Having in mind further applications of it  we  shall
prove it in slightly more general form than is needed now.)

{\bf6.30. Proposition.} {\it Let $C$ be a convex set in a Banach space $X$ and let
$y^{**}\in w^*-cl(C)\subset X^{**}$ and $||y^{**}-x||<1$
for some $x\in $X.
Then there  exists $y\in C$ such that}  $||y^{**}-y||<2.$

Proof. Let  real  number $\delta $  be  such  that $||y^{**}-x||\le
\delta <1,$  i.e.
$y^{**}-x\in \delta B(X^{**})$.  Then $x=(x-y^{**})+y^{**}\in
\delta B(X^{**})+w^*-cl(C)\subset w^*-cl(\delta B(X)+C)$.
Hence $x\in w-cl(\delta B(X)+C)$. Since $\delta B(X)+C$ is convex, 
then  we can  find
$u\in \delta B(X)$ and $y\in C$ such that $||x-(u+y)||<1-\delta $.
Then $||y-x||<||u||+1-\delta \le 1.$ Hence
$||y^{**}-y||\le||y^{**}-x||+||x-y||<2.$ Proposition is proved.

{\bf6.31.} Let $X$ be arbitrary Banach space and $Y\in G(X)$.
Let us denote by $Q$ the quotient mapping $Q:X^{**}\to X^{**}/X.$
Applying Proposition 6.30  to $C=Y$
we obtain the following statement.

{\bf Corollary.} {\it The subspace $Q(Y^{\perp \perp })$  belongs
to $G(X^{**}/X)$  and  is isomorphic to} $Y^{**}/Y$.

{\bf6.32.} Each of the  statements  of  Theorem  6.29  follows  from  the
comparison of corollary  6.31  and  the  following  lemma.  (We  use
notation of 6.31).

{\bf Lemma.} {\it For arbitrary Banach space $X$  and  arbitrary  subspaces
$Y,Z\in G(X)$ we have} $\Theta_0(Q(Z^{\perp \perp }),Q(Y^{\perp \perp }))\le
2\Theta_0(Z,Y)$.

Proof. Let $\varepsilon >0$ be arbitrary. Let
$\xi \in S(Q(Z^{\perp \perp }))$. Then we can  find
$z^{**}\in Z^{\perp \perp }$ and $x\in X$ such
that $Qz^{**}=\xi $  and $||z^{**}-x||<1+\varepsilon $.  By  Proposition
6.30 we can find $z\in Z$ such that $||z^{**}-z||<2(1+\varepsilon )$.
Since $\Theta_0(Z^{\perp \perp },Y^{\perp \perp })=
\Theta_0(Z,Y)$ (by Theorem 3.4(d)), then we  can  find
$y^{**}\in Y^{\perp \perp }$  such  that
$||z^{**}-z-y^{**}||<2(1+\varepsilon )\Theta_0(Z,Y)$.
Hence $||\xi -Qy^{**}||<2(1+\varepsilon )\Theta_0(Z,Y)$.
Since $\varepsilon >0$
is arbitrary, we obtain the desired inequality. Lemma is proved.

{\bf6.33. Definition.} Let $P$ be a property of Banach spaces.
Class $Q$  is said to be the {\it preproperty} of $P$ (Notation:
$Q$=pre($P$)) if
$$
(X\in Q)\Leftrightarrow (X^*\in P).
$$

{\bf6.34. Proposition.} {\it If $P$ is stable property then {\rm pre}$P$
is also stable. If $P$ is extendedly stable property then {\rm pre}$P$
is also  extendedly  stable.  Moreover  we  have
$s(${\rm pre}($P))\ge s(P)$ and} $es(${\rm pre}($P))\ge es(P)$.

Proof. Let $X$ be a Banach space and let $Y,Z\in G(X)$  be  such  that
$Y\in $pre($P$) and $\Theta_0(Z,Y)<es(P)$. Let us prove that $Z\in $pre($P$).
In order  to
do that let us consider space $X_1=X\oplus _{2}Y$  and
let $Y_{t} (t>0)$  be  the
following family of subspaces:
$$
Y_{t}=\{(y,ty):y\in Y\}.
$$
It is clear that every $Y_{t}$ is isomorphic (even isometric)  to $Y$
and $\lim_{t\to 0}\Theta (Y,Y_{t})=0.$ Therefore choosing $t$  small  enough
we  obtain
$\Theta_0(Z,Y_{t})<es(P)$. On the other hand we have $Z\cap Y_{t}=0$
and $\delta (Z,Y_{t})>0.$  Let
$X_{2}={\rm lin}(Z\cup Y_{t})$. By remark from 5.3 $X_{2}$ is a closed
subspace of $X$  and
$Y_{t}$ is a complemented subspace of it. So we have:

(a) $Y^{\perp }_{t}$ is isomorphic to $Z^*$.

(b) $Z^{\perp }$ is isomorphic to $(Y_{t})^*$.

By duality formula (Theorem 3.4(d)) we have

(c) $\Theta_0(Y^{\perp }_{t},Z^{\perp })<es(P).$

Since $(Y_{t})^*\in P$ and $P$ is isomorphic invariant by Proposition 6.7,
we obtain that $Z^*\in P$. So $Z\in $pre($P$). Proposition is proved.

{\bf6.35. Remark.} {\it Analogous statement can be proved for open
properties
which are  isomorphic  invariants.  But  the  known  proof  of  this
statement is rather long (see 7.12).}

{\bf6.36.} The following conjecture seems to be excessively audacious but
at the moment I don't no counterexamples.

{\bf Conjecture.} {\it For every cardinal $\alpha $  the
intersection  of  every
extendedly stable property  with  the  set  of  Banach  spaces  with
density character less than $\alpha $ is an $l_1$-property.}

We need to consider intersections here in order to  avoid  such
trivial counterexamples as the class of all Banach spaces.

It  should  be  noted  that  stable  properties  need  not   be
$l_1$-properties. (It follows because the complement of stable property
is stable and the complement of $l_1$-property is not an $l_1$-property.)

{\bf6.37.} In [O9] the present author obtained several results in support
of this conjecture. In this connection the following  definition  is
useful.

{\bf Definition.} Class $P$  of  Banach  spaces  is  called  a  {\it regular}
$l_1${\it -property} if there exist a real
number $\delta >0,$ set $\Gamma $  and  a  subset
$A\subset S(l_1(\Gamma ))$ satisfying the conditions:

 1. The set $A$ consists of finitely non-zero vectors.

 2. If $a_0\in A$, then $A$ contains all vectors $a\in S(l_1(\Gamma ))$
for which
$$
(\forall \gamma \in \Gamma )({\rm sign}a_0(\gamma )={\rm sign}a(\gamma )).
$$

3. For a Banach space $X$ the following conditions are equivalent:

(a) $X\not\in P.$

(b) $h(X,A)>0.$

(c) $h(X,A)\ge \delta .$

Supremum of those $\delta >0$ for which there exist $\Gamma $
and $A\subset S(l_1(\Gamma ))$  such
that conditions 1--3 are satisfied is called the {\it regular exponent}  of
$P$.

{\bf6.38.} In [O9] the following results were proved.

(a) Properties listed in 6.27 are regular $l_1$-properties.

(b) The union of every set of regular $l_1$-properties  with  uniformly
bounded away from zero regular exponents is a regular $l_1$-property.

(c) If $P$ is regular $l_1$-property then  pre($P$)  and $P^{co}$  are  regular
$l_1$-properties.

I do not know examples of $l_1$-properties which are not  regular.
In particular I do not know  whether  any  intersection  of  regular
$l_1$-properties  with  uniformly  bounded  away  from   zero   regular
$l_1$-exponents is a regular $l_1$-property.

\medskip
\centerline{{\bf Notes and remarks}}

{\bf6.39.} Let us give some historical comments.

 Proposition 6.1 and the part of Proposition 6.2 which  concerns  the
operator opening seems to be new. The second part of proposition 6.2
is taken from [O1]. It was also announced in [Fr].  Observation  6.3
is due to A.Douady [D]. Construction of 6.14  is  a  straightforward
generalization of M.I.Kadets' construction [Kad2]. Approach of  6.17
was introduced by the author [O2], [O3], [O4]. Proposition  6.21  is
taken from [O4] (see also [O3]).

Corollary 6.22 is a slight strengthening of the  following
result due to A.L.Brown [Br1]: Class of  reflexive  spaces  is  stable
with the stability exponent no less than  $\sqrt2-1.$

Corollary 6.23 was proved in [KKM].

The main idea of the description of  super-reflexivity,  which
was  mentioned  in  6.27.6  is  due  to  R.C.James [J2].  Nessesary
additions was made in [J3], [SS] and [JS].  Many  other  result  on
super-reflexive spaces can be found in [Be2] and [Du].

Papers [AAG] and [Ja2] contains  somewhat  weaker  versions  of
Proposition 6.34 with more complicated proofs. The present proof  is
due to the author.

It turns out that the paper [O6] which is relevant to the topic
of this chapter is ``almost free'' of  new  results  and  is  only  of
historical interest. This paper contains proofs of Proposition  6.11
and of stability of quasireflexivity.

{\bf6.40.} I would like to add the following information  concerning  the
statement of 6.28 about the distinction between $l_1$-properties.  The
direct sum $(\sum^{\infty }_{n=1}\oplus l^n_1)_{2}$ is an example
of reflexive non-$B$-convex space.

Examples  of  non-reflexive   and   hence   non-super-reflexive
$B$-convex spaces was constructed by R.C.James [J4], see also [JL].

It is known [BrS] that $c_0\in ABS$. Hence $ABS$ does not imply any  of
the properties: reflexivity, $B$-convexity, super-reflexivity, $BS$.

A.Baernstein [Bae] constructed reflexive space without $BS$.

{\bf6.41.} In connection with definition 6.18 and Proposition 6.21 it  is
natural to propose the problem of description of the sets
$$\sigma (A)=\{h(X,A):\ X\hbox{ running through the class of all Banach
spaces}\}$$
for different sets $A\subset S(l_1(\Gamma ))$ and to calculate $h(X,A)$
for classical Banach spaces.

In addition to results mentioned in 6.19, 6.20 and 6.27 I  know
results of the mentioned type only for the set $A\subset S(l_1)$
described  in
6.27.3 and its uncountable analogues  mentioned  in  6.27.2.  It  is
known that $h(l_{p},A)=2^{1/p}/2$ for $1\le p<\infty$  [BRR] (see
also [Ko1] and [WW], p.~91). The  book [WW]  contains
also  results  on  evaluation  of
$h(L_{p}(\mu ),A)$, where $\mu $ is arbitrary measure.

In addition, J.Elton and E.Odell [EO] (see also [Di,  p.~241])
proved that $h(X,A)>1/2$ for every infinite dimensional  Banach  space
$X$. Hence in this case $\sigma (A)=\{0\}\cup (1/2,1]$.

For uncountable analogues of $A$ (see 6.27.2) there exist  spaces
for which the value of index $h$ is 1/2 (see [EO]). So  in  this  case
the set of possible values of index $h$ is $\{0\}\cup [1/2,1]$.

Some  results  closely  connected  with  mentioned  above   are
contained in [Ba] and [Ko2].

 {\bf6.42.} M.~G.~Krein, M.~A.~Krasnoselskii  and  D.~P.~Milman  [KKM,  p.~104]
asked the following question: does inequality $\Theta (Y,Z)<1$  imply  that
dens$Y=$dens$Z$? This question was solved in negative in [O1],  where  a
Banach space $X$ was  constructed  such  that  for  certain  subspaces
$Y,Z\in G(X)$ we have $\Lambda (Y,Z)\le 2\sqrt2-2,\ Z$ is separable and
dens$Y$  equals  to the power of continuum. Here we reproduce this example.
Let  $a=\sqrt2-1$. Consider the algebraic sum
$$
X=c_0([0,1])\oplus (\sum^{\infty }_{i=1}
\oplus C(0,1))_1\oplus (\sum^{\infty }_{j=1}\oplus C(0,1))_1.
$$

We endow it with the norm:
$$
||(h_0,(h_{i})^{\infty }_{i=1},(g_{j})^{\infty }_{j=1})||=
$$
$$
\max (\sum^{\infty }_{i=0}||h_{i}-(1/2)g_{i+1}||,\
\sum^{\infty }_{j=1}||g_{j}-(1/2)h_{j}||),
$$
where all the norms on  the  right-hand  side  are  suprema  of  the
modulus on [0,1]. It is clear that the space
$$Y=\{(h_0,(h_{i})^{\infty }_{i=1},(0)^{\infty }_{j=1})\}$$
is   isometric   to $c_0([0,1])\oplus_1(\sum^{\infty }_{i=1}\oplus C(0,1))_1$
and    the    space
$$Z=\{(0,(0)^{\infty }_{i=1},(g_{j})^{\infty }_{j=1})\}$$
is isometric
to $(\sum^{\infty }_{j=1}\oplus C(0,1))_1$.  It  is  also
clear that $Z$  is  separable  and  dens$Y$  equals  to  the  power  of
continuum. Let us  estimate $\Lambda (Y,Z)$.  To  do  this  it  suffices  to
estimate dist$(y_{i},B(Z)) (i\in {\bf N}\cup \{0\})$ and
dist$(z_{j},B(Y)) (j\in {\bf N})$ for vectors
$y_{i}\in S(Y)$ of the form
$$
y_0=(h_0,(0)^{\infty }_{i=1},(0)^{\infty }_{j=1})\hbox{ and }
y_{i}=(0,(0,\ldots ,0,h_{i},0,\ldots ),(0)^{\infty }_{j=1}),
$$
where $h_{i}$ is at the $i$-th position, and
$$
z_{j}=(0,(0)^{\infty }_{i=1},(0,\ldots ,0,g_{j},0,\ldots )),
$$
where $g_{j}$ is at the $j$-th position. This is done in  the  same  manner
for all such vectors except for vectors of the form $y_0$.

We verify the estimate for $y_1$. Consider the vector
$$
f=(0,(0)^{\infty }_{i=1},(0,ah_1,0,\ldots ))\in B(Z).
$$
We have
$$
\hbox{dist}(y_1,B(Z))\le||y_1+f||=\max \{||(1-a^{2})h_1||,||ah_1||+||ah_1||\}=
$$
$$
\max \{(1-a^{2}),2a\}=2\sqrt2-2.
$$
Let us now verify the estimate for $y_0$. Let $\varepsilon >0.$
We  introduce
the set $A=\{x$: $h_0(x)$ $\ge \varepsilon \}$. This set is finite.
Let the function $w_{\varepsilon }:A\to {\bf R}$
be defined as $w_{\varepsilon }(x)=h_0(x)$ and extend it as a continuous
function  to
$[0,1]$ in such a way  that
$\sup \{$ $w_{\varepsilon }(x)$ :$x\in [0,1]\}=1.$  Then  for  every
$b\in [0,1]$ we have
$$
||h_0-bw_{\varepsilon }||=\max \{1-b,b+\varepsilon \}
\eqno{(12)}$$
Let $z=(0,(0)^{\infty }_{i=1},(aw_a,0,\ldots)).$ We have $||z||\le a$ and
dist$(y_0,B(Z))\le||y_0+z||=\max \{||h_0-a^{2}w_{a}||,||aw_{a}||\}=
\max \{1-a^{2},a^{2}+a,a\}=2\sqrt2-2.$

{\bf6.43. Remark.} {\it Inequality   (12)   implies   that   for   subspaces
$c_0([0,1])\subset l_{\infty }([0,1])$ and $C(0,1)\subset l_{\infty }([0,1])$
we have}
$$
\Lambda _0(c_0([0,1]),C(0,1))=1/2.
\eqno{(13)}$$

Therefore the estimates of extended stability  exponents  in  6.27.2
are precise when $\alpha $ is not greater then the cardinality of continuum.

Equality (13) was independently  noticed  by  A.N.Plichko,  who
used it in the theory of biorthogonal systems [Pl].

{\bf6.44.} The author proved [O1] that if $Y$ and $Z$  are  subspaces  of  a
Banach space with an extended unconditionally  monotone  basis  (see
definition in [S2], \S 17), then $\Theta (Y,Z)<1$ implies dens$Y$=dens$Z$.

V.I.Gurarii [Gu2] proved that if we suppose that $X$ is uniformly
convex (=uniformly rotund) and $\delta $ is its modulus  of  convexity  (see
necessary definitions in [Da]  or  [Di]), 
then $\Theta (Y,Z)<1/2+\delta (1/7)/2$
implies dens$Y$=dens$Z$.

 {\bf6.45.} The examples from 6.42 has predecessor (see [Le]  and [S3,
p.~271]).

{\bf6.46.} Starting from [KKM] the following dual version of the  problem
to which the present section is devoted  is  considered:  How  close
should be the structure of the quotients over  the  subspaces  which
are close with respect to the geometric opening?

{\bf Definition.} A property $P$ is called {\it co-open} if for every  Banach
space $X$ the set of those $Y\in G(X)$ for  which $X/Y\in P$  is  open  in  the
topology induced by $\Theta $.

A property $P$ is called {\it costable} if there exists  a  number $\alpha >0$
such that for every Banach space $X$ and every $Y,Z\in G(X)$, if $X/Y\in P$  and
$\Theta (Y,Z)<\alpha $, then $X/Z\in P$. The least upper bound of
numbers $\alpha $  for  which
this statement is true is denoted by $c(P)$.

A property $P$ is called {\it extendedly costable} if  there  exists  a
number $\alpha >0$ such that for every Banach space $X$ and every
$Y,Z\in G(X)$, if
$X/Y\in P$ and $\Theta_0(Y,Z)<\alpha $, then $X/Z\in P$. The least upper
bound of numbers $\alpha $
for which this statement is true is denoted by $ec(P)$.

It turns out that the concept of costability coincides with the
concept  of  stability  and  the  concept  of  extended  costability
coincides with the concept of extended stability. Moreover, we  have
$s(P)=c(P)$ and $es(P)=ec(P)$ for every property $P$. This  assertion  can
be proved using arguments of Proposition 6.34 and  ``co-analogue''  of
Proposition 6.2. Somewhat weakened version  of  this  statement  was
proved by J.~Alvarez, T.~Alvarez  and  M.~Gonzalez  [AAG].  Some  other
predecessors of this result can be found in [O2].

If  we  restrict  ourselves  by  consideration  of   isomorphic
invariants then the concept of open property would coincide with the
concept of co-open property. The  known  proof  of  this  result  is
rather complicated (see [Ja2]).

 {\bf6.47.} Proposition 6.11 (a) implies that if property $P$ is  (a)  open;
(b) isomorphic invariant; (c) closed with respect  to  formation  of
direct sums; then the solution of the three space problem for $P$  is
positive. In particular it is so for  stable  properties  which  are
closed with respect to formation of direct sums. Here it  should  be
noted that the three space problems for properties  listed  in  6.27
(with obvious exception of 6.27.4) were already solved. Here is  the
list  of  corresponding  references.  Reflexivity: [KS,  p.~575];
$B$-convexity [Gi, p. 103] (see also [ELP], [R]);  super-reflexivity:
[ELP] (see also [H, p.~86], [R], [Sc]), class of  spaces,  which  do
not  contain  isomorphic  copies  of $l_1$: [R] (see  also   [AGO]),
alternate-signs Banach-Saks property: [O3];  Banach-Saks  property:
[GR] (see also [O3]).

 6.48. The papers [Ja2], [O1], [O3] (see also [O6])  deals  with  the
following question: for what Banach  spaces $Z$  classes $P_{Z}=\{X:X$  is
isomorphic to $Z\}$ are open and what is the least upper bound  of $\alpha >0$
for which the assertion of Proposition 6.6 is true if  we  let $P=P_{Z}$
and $Y=Z$? The most interesting problem seems to be the following:  is
it true that every $Z$ for which $P_{Z}$ is open  is  either  injective  or
isomorphic to $c_0
$ or $l_1(\Gamma )$? The fact  that $P_{Z}$  is  open  when $Z$  is
injective can be easily derived from the Proposition 5.3.

 {\bf6.49.} There is another survey devoted to the same topic as section 6
of the present survey.  It  is  the  survey  due  to  J.~Alvarez  and
T.~Alvarez [AA].  Unfortunately  because  of  certain  confusion  in
different concepts of stability,  some  statements  of  this  survey
(Theorem 12, Observation 15) are incorrect.

{\bf6.50.} I am sure that the results and the methods  of  papers  [DJL],
[DL], [Fa], [J5] and [PX]
 will be useful in  further  investigations
of $l_1$-properties.

{\bf7. Additional remarks on the topologies on the set of all  subspaces}

{\bf of a Banach space and their applications}

{\bf7.1.} Basic facts about geometric opening are discussed in well-known
course of I.M.Glazman and Yu.I.Lyubich [GL].

{\bf7.2.} Openings found applications in the theory of best approximation
in Banach spaces. Recall that a subspace $Y$ of a Banach  space $X$  is
called a {\it Chebyshev subspace} if for every $x\in X$ there exists
unique $y\in Y$
such that $||x-y||$=dist$(x,Y)$. Subspace $Y$ is called {\it almost-Chebyshev}
if the condition above is satisfied for all $x$ except some  set  of  the
first category. It is known [Ga] that $c_0$ does  not  contain  almost-
Chebyshev proper infinite dimensional  subspaces.  A.L.Garkavi [Ga]
proved that if a separable  Banach  space $X$  contains  a  reflexive
infinite dimensional subspace then it contains  an  almost-Chebyshev
proper  infinite  dimensional  subspace.  He  also  proved  that   a
separable  dual  Banach  space  contains   almost-Chebyshev   proper
infinite dimensional subspace. The proof  relies  heavily  on  Baire
category argument for metric space $(G(X),d_{op})$. Using Corollary 6.22,
A.L.Brown [Br1] showed that the first result  can  be  proved  using
category argument for the metric space $(G(X),d_g)$.

Some other applications of  openings  in  the  theory  of  best
approximation can be found in [S3].

{\bf7.3.} V.P.Fonf [Fo] used geometric opening in his  investigations  of
polyhedral Banach spaces. Recall  that  a  Banach  space  is  called
{\it polyhedral} if the unit  ball  of  each  of  its  finite  dimensional
subspaces is a polytope. V.P.Fonf [Fo] proved  that  every  infinite
dimensional polyhedral Banach space contains an isomorphic  copy  of
$c_0$.

{\bf7.4.} The topology induced by the opening on the set of all subspaces
of a Hilbert space was used by R.G.Douglas and  C.M.Pearcy [DP]  in
investigations of the lattice  of  invariant  subspaces  of  bounded
linear operators in Hilbert space.

{\bf7.5.} Theorem 3.4 (e)  has  many  applications  in  constructions  of
``nicely bounded'' biorthogonal systems  (see [Da, p.~93],  [Pel3],
[Te]).

{\bf7.6.} Geometric opening is a natural tool in some  questions  of  the
theory of bases (see [S1], [S2]).

{\bf7.7.} The topology on $G(X)$ induced by the geometric opening  appeared
in a natural way in investigations of infinite dimensional  analytic
manifolds [D], operator function equations [Ja1], [KT], [Man], [Sl].

{\bf7.8.} The  geometric  opening  and  quantity $\Theta _0$  turned  out
to  be
important tool in the creation of the theory of  Fredholm  operators
in Banach spaces [GK], [K1]. They are repeatedly used in the  theory
of operators in Banach spaces  (see,  in  particular, [AZ], [CPY],
[GM2], [Go], [GoM], [MeS], [Sob], [Va2], [Va3]), Fredholm and  semi-Fredholm
complexes of Banach spaces (see [AV], [A], [Do], [F], [FS],
[Va1]).

{\bf7.9.} Let $X$ and $Y$ be Banach spaces. If we have a topology
on $G(X\oplus Y)$
then, considering graphs we obtain a topology  on  the  set  of  all
closed linear operators with domains  in $X$  into  $Y$.  The  topology
obtained in such a way from the topology induced  by  the  geometric
opening turned out to be very useful. Many important characteristics
of closed linear operators are stable with respect to this topology.
This circle of problems was  investigated  by  J.D.Newburgh  [New1],
[New2], H.O.Cordes and J.-P.Labrousse [CL], G.Neubauer [N1], [N2],
T.Kato [K2], J.-P.Labrousse [L1], [L2]. Some of  these  results  are
presented in the book due to K.R.Partasarathy [Par]. If we  restrict
ourselves by consideration of self-conjugate operators in a  Hilbert
space then this topology coincides with the topology of uniform
resolvent convergence (see [RS], chapter VIII).

It should be mentioned that some ``unbounded'' analogues of known
results fails for this topology (see [L3] and forthcoming  book  due
to J.-P.Labrousse). In the case of unbounded closed densely  defined
operators in Hilbert space the situation sometimes can be  saved  by
the use of another metric (see [L3], [LM] and forthcoming  book  due
to J.-P.Labrousse). This metric goes back to  C.Davis  [Dav].  (This
metric is defined only on the subset of $G(H\oplus H)$ consisting of  graphs
of closed densely defined operators on the separable  Hilbert  space
$H$.)

{\bf7.10.} Let $X,Y$ be Banach spaces and $T\in L(X,Y)$. The
{\it minimum modulus}  of
$T$ is defined by
$$
\gamma (T)=\sup \{c\ge 0:(\forall x\in X)
(||Tx||\ge c\cdot \hbox{dist}(x,\hbox{ker}T))\}.
$$
 (We refer to [K2] for basic properties of $\gamma )$.

A.S.Markus [M1] proved the following result.

{\bf Theorem.} {\it Let $T,S\in L(X,Y)$. Then}

(a) $\Theta _0(\hbox{ker}S,\hbox{ker}T)\le \gamma (T)^{-1}||S-T||.$

(b) $ \Theta _0(\hbox{im}S,\hbox{im}T)\le \gamma (S)^{-1}||S-T||.$

(c) {\it If} $\Theta ($ker$S$,ker$T)<1/2,$ {\it then}
$$
|\gamma (S)-\gamma (T)|\le 3||S-T||(1-2\Theta (\hbox{ker}S,\hbox{ker}T))^{-1}.
$$
(d) {\it If} $\Theta ($im$S$,im$T)<1/2,$ {\it then}
$$
|\gamma (S)-\gamma (T)|\le 3||S-T||(1-2\Theta (\hbox{im}S,\hbox{im}T))^{-1}.
$$

{\bf7.11.} R.Janz [Ja2] proved the following statement. Let $D$ be a metric
space and $T:D\to G(X)$ be a mapping continuous in the  topology  induced
by the geometric opening. Then there exist Banach spaces $Z_1$  and $Z_2$
and continuous mappings $R_1:D\to L(Z_1,X)$ and $R:D\to L(X,Z_2)$ such  that
for every $d\in D$ we have

(a) ${\rm im}(R_1(d))=T(d).$

(b) $ \ker (R_2(d))=T(d).$

(c) $||R_1(d)||\le 1,\ ||R_2(d)||\le 1.$

(d) $ \gamma (R_1)\ge 1/10, \gamma (R_2)\ge 1/10.$

{\bf7.12.} Applying 7.11 to the identity mapping $T:G(X)\to G(X)$  and  using
7.10 it is not hard to prove assertions formulated in 6.35 and 6.46.

{\bf7.13.} Let $X$ be a linear space with  two  norms  and $Y,Z$  be  linear
subspaces of $X$. By $\Theta _1(Y,Z)$ and $\Theta _2(Y,Z)$ we shall denote
the geometric
openings  in  the  sense  of  the  first  and   the   second   norms
respectively. (Here we use the fact that the  geometric  opening  is
well-defined for non-closed subspaces of non-complete normed spaces.

A.S.Markus [M2] proved the following statement: if for some $c>0$
and every linear subspaces $Y,Z\subset X$ we have
$$\Theta _1(Y,Z)\ge c\Theta _2(Y,Z),$$
then the norms are equivalent.

{\bf7.14.} Some properties of the metric space $(G(H),d_g)$, where $H$ is  the
Hilbert space were considered in [I1], [I2], [Lo1], [Lo2], [Rod].

{\bf7.15.} H.Porta and L.Recht [PR] proved that for every Banach space $X$
there exists a mapping $\phi :G_c(X)\to G_c(X)$ such that for every
$Y\in G(X)$  the
space $\phi (Y)$ is a complement of $Y$ and $\phi $
is continuous with respect  to
the $d_{op}$.

\centerline{REFERENCES}

[AG] N.I.Akhiezer and I.M.Glazman, Theory of linear  operators  in
Hilbert space, Moscow, Gostekhizdat, 1950 (in Russian); see also
English translation of the third edition of  this  book:  Pitman
Press, 1980.

[AV] E.Albrecht and F.-H. Vasilescu, Stability of  the  index  of  a
semi-Fredholm complex of Banach spaces J. Funct. Anal. 66  (1986),
141--172.

[AA]  J.Alvarez  and   T.Alvarez,   Estabilidad   de   propiedades
isomorphas  de  espacios  de  Banach  respecto  del  gap,  Univ.
Cantabria, Serv. de Publ, Santander, 1991, pp. 1--8.

[AAG]  J.Alvarez,  T.Alvarez  and  M.Gonzalez,  The  gap   between
subspaces and perturbation of non semi-Fredholm operators, Bull.
Austral. Math. Soc., 45 (1992), no. 3, 369--376.

[AGO] T.Alvarez, M.Gonzalez and V.M.Onieva, A note on  three-space
Banach space ideals, Arch. Math. 46 (1986), 169--170.

[A] C.-G.Ambrozie, Stability of the index of a Fredholm  symmetrical
pair, J. Operator Theory, 25 (1991), no. 1, 61--77.

[AZ] B.Aupetit and J.Zemanek, Uniformly regular families  of  commuting
operators J. Funct. Anal. 78 (1988), 24--30.

[Bae] A.Baernstein, On reflexivity and summability, Studia Math.  42
(1972), 91--94.

[Ba] K.Ball, Inequalities  and  sphere-packing  in $l_p$ ,  Israel  J.
Math., 58 (1987), no. 2, 243--256.

[Be1] B.Beauzamy,  Banach-Saks  properties  and  spreading  models,
Math. Scand., 44 (1979), no. 2, 357--384.

[Be2]  B.Beauzamy,  Introduction  to  Banach  Spaces   and   their
Geometry, Amsterdam: North-Holland Publishing Company, 1982.

[BL] B.Beauzamy et J.-T. Lapreste, Modeles etales des  espaces  de
Banach, Paris, Hermann, 1983.

[BDGJN] G.Bennett, L.E.Dor, V.Goodman, W.B.Johnson and C.M.Newman,
On uncomplemented subspaces of $L_p , 1<p<2,$  Israel  J.Math.,  26
(1977), no. 2, 178--187.

[Ber] E.Berkson, Some metrics on the subspaces of a Banach  space,
Pacif. J. Math., 13 (1963), no. 1, 7--22.

[Bor] K.Borsuk, Drei  Satze  uber  die $n$-dimensionale  euklidishe
Sphare, Fund. Math. 20 (1933), 177--190.

[B]  J.~Bourgain, A  counterexample   to   a   complementation   problem,
Compositio  Math., 43 (1981), no. 1, 133--144.

[Br1] A.L.Brown, On the space of subspaces of a Banach  space, J.
London Math. Soc., 5 (1972), no. 1, 67--73.

[Br2] A.L.Brown,  The  Borsuk-Ulam  theorem  and  orthogonality  in
normed linear  spaces,  Amer.  Math.  Mon.,  86  (1979),  no.  9,
766--767.

[BrS] A.Brunel and L.Sucheston, On $J$-convexity and  some  ergodic
super-properties of Banach spaces, Trans.  Amer.  Math.  Soc.  204
(1975), 79--90.

[BRR] J.A.C.Burlak, R.A.Rankin and A.P.Robertson, The  packing  of
spheres in the space $l_p$ , Proc. Glasgow  Math. Assoc. 4 (1958),
22--25.

[CPY] S.R.Caradus, W.E.Pfaffenberger and B.Yood,  Calkin  algebras
and algebras of operators on Banach  spaces,  New  York,  Marcel
Dekker, 1974.

[CG1] J.M.F.Castillo and M.Gonzalez, The  Dunford-Pettis  property
is not a three-space property, preprint.

[CG2] J.M.F.Castillo and M.Gonzalez, Properties ($u$)  and  ($V$)  are
not three-space, preprint.

[CL] H.O.Cordes and J.P.Labrousse, The invariance of  the  index  in
the metric space of closed operators, J. Math. Mech.,  12  (1963),
693--719.

[C] H.Corson, The weak topology of  a  Banach  space,  Trans.  Amer.
Math. Soc., 101 (1961), no. 1, 1--15.

[Dav] C.Davis, Separation of two linear subspaces, Acta sci. math.,
19 (1958), no. 3--4, 172--187.

[DJL] W.J.Davis, W.B.Johnson and J.Lindenstrauss, The $l^n_1$-problem and
degrees of non-reflexivity, Studia Math. 55 (1976), 123--139.

[DL] W.J.Davis and J.Lindenstrauss, The $l^n_1$-problem  and  degrees  of
non-reflexivity, II, Studia Math. 58 (1976), 179--196.

[Da] M.M.Day, Normed linear spaces, Berlin, Springer-Verlag, 1973.

[DKL] R.A.DeVore, H.Kierstead and G.G.Lorentz, A proof of Borsuk's
theorem, Lect. Notes Math., 1332 (1988), 195--202.

[Di] J.Diestel, Sequences and series in Banach spaces,  New  York,
Springer-Verlag, 1984.

[D] A.Douady, Le probleme des modules  pour  les  sous-espaces  d'un
espace analytique donne, Ann. Inst. Fourier,  15  (1966),  no.  1,
1--94.

[Do] R.Douady, Petites perturbations d'une  suite  exacte  et  d'une
suite quasi-exacte, Seminaire d'Analyse, Nice, 1965--66, pp. 21--34.

[DP]  R.G.Douglas  and  C.M.Pearcy,  On  a  topology  for  invariant
subspaces, J. Funct. Anal., 2 (1968), 323--341.

[Du] D.van Dulst,  Reflexive  and  superreflexive  Banach  spaces,
Amsterdam, Mathematish Centrum, 1978.

[EO] J.Elton  and  E.Odell,  The  unit  ball  of  every  infinite-
dimensional  normed  linear  space  contains a $(1+\varepsilon )$-separated
sequence, Colloquium Math., 44 (1981), no. 1, 105--109.

[ELP] P. Enflo, J.Lindenstrauss and G.Pisier, On the  ``three-space
problem'', Math. Scand. 36 (1975), 199--210.

[F] A.~S.~Fainshtein, On the opening
between
subspaces
of  a  Banach
space, manu- script, Baku, 1983 (in Russian).

[FS] A.S.Fainshtein and V.S.Shul'man, Stability of  the  index  of
short  Fredholm  complex  of  Banach  spaces  with  respect   to
perturbation, small with respect to measure of  non-compactness,
Spectral Theory of  Operators  (Baku),  4  (1982),  189--198  (in
Russian).

[Fa] J.E.Farahat, On the problem of $k$-structure, Israel J. Math.  28
(1977), no. 1--2, 141--150.

[Fo] V.P.Fonf, On  polyhedral  Banach  spaces,  Mat.  Zametki,  30
(1981), no. 4, 627--634 (in Russian).

[Fr] F.J.Freniche, Una distancia entre espacios de Banach, Actas IX
Jornadas  Matematicas  Hispano-Lusas,  Vol.  I,  Universidad   de
Salamanca, 1982, pp. 281--283.

[Ga] A.L.Garkavi, On  Chebyshev  and  almost-Chebyshev  subspaces,
Izv. Akad. Nauk SSSR, Ser. Mat. 28 (1964), no.  4,  799--818  (in
Russian).

[Gi] D.P.Giesy, On a convexity condition in normed linear  spaces,
Trans. Amer. Math. Soc. 125 (1966), no. 1, 114--146.

[GL] I.M.Glazman, Yu.I.Lyubich, Finite dimensional linear  analysis,
Moscow, Nauka, 1969 (in Russian).

[GR] B.V.Godun and S.A.Rakov, Banach-Saks property and  the  three
space problem, Mat. Zametki, 31 (1982), 61--74 (in Russian).

[GK] I.C.Gohberg and M.G.Krein, The basic propositions  on  defect
numbers, root numbers and indices of linear  operators,  Uspekhi
Mat. Nauk, 12  (1957),  no.  2,  43--118  (in  Russian);  English
translation: Transl. Amer. Math. Soc. (2), 13 (1960), 185--264.

[GM1] I.C.Gohberg and A.S.Markus,  Two  theorems  on  the  opening
between subspaces of a  Banach  space,  Uspekhi  Mat.  Nauk,  14
(1959), no. 5, 135--140 (in Russian).

[GM2] I.C.Gohberg and  A.S.Markus,  Characteristic  properties  of
certain points of the spectrum of bounded linear operators, Izv.
VUZ, Ser. Mat. (1960), no. 2, 74--87 (in Russian).

[Go] M.Gonzalez, A  perturbation  result  for  generalized  Fredholm
operators in the boundary of the group  of  invertible  operators,
Proc. R. Ir. Acad. 86A (1986), no. 2, 123--126.

[GoM] M.Gonzalez and A.Martinon, Operational quantities derived from
the norm and generalized Fredholm  theory,  Comment.  Math.  Univ.
Carolinae 32 (1991), no. 4, 645--657.

[Gu1] V.I.Gurarii, Openings and inclinations  of  subspaces  of  a
Banach space, Teor Funktsii, Funktsional. Anal. i  Prilozhen.  1
(1965), 194--204 (in Russian).

[Gu2] V.I.Gurarii, On uniformly convex and uniformly smooth Banach
spaces,  Teor.  Funktsii,  Funktsional.  Anal.  i  Prilozhen.  1
(1965), 205--211 (in Russian).

[GuM] V.I.Gurarii and A.S.Markus, On the  geometric  and  operator
definitions of the opening between  subspaces,  Issled.  Algebr.
Mat.
Anal.,
Kishinev,
Kartya
Moldovenya- ske, 1965,  pp.  128--130.
Correction: Mat. Issled 2 (1967), no. 1, 169 (in Russian).

[H] S.Heinrich, Ultraproducts in Banach space theory, J. reine und
angew. Math., 313 (1980), 72--104.

[I1] S.Izumino, Inequalities on transforms of subspaces of a Hilbert
space, Math. Japon., 25 (1980), no. 1, 131--134.

[I2] S.Izumino,  Generalized  inverses  method  for  subspace  maps,
Tohoku Math. J. (2) 35 (1983), no. 4, 649--659.

[J1] R.C.James, Weak compactness and reflexivity, Israel J. Math.,
2 (1964), 101--119.

[J2] R.C.James, Uniformly non-square Banach spaces, Ann.  Math. 80
(1964), no. 3, 542--550.

[J3] R.C.James, Some self-dual properties of normed linear spaces,
Sympos. Infinite-Dimensional Topology (Baton  Rouge, La, 1967)
Ann.  of  Math.  Studies,  vol.  69,  Princeton   Univ.   Press,
Princeton, N.J., 1972, pp. 159--175.

[J4] R.C.James A nonreflexive Banach space that  is  uniformly
non-octahedral, Israel J. Math. 18 (1974), 145--155.

[J5] R.C.James, A separable somewhat  reflexive  Banach  space  with
nonseparable dual, Bull.  Amer.  Math.  Soc.  80  (1974),  no.  4,
738--743.

[JL] R.~C.~James  and  J.~Lindenstrauss,  The  octahedral  problem  for
Banach spaces, in:``Proceedings of the seminar on  random  series,
convex sets and geometry of Banach spaces'', Aarhus, Denmark, 1974,
pp. 100--120.

[JS] R.C.James and J.J.Sch\"affer, Super-reflexivity and  the  girth
of spheres, Israel J. Math., 11 (1972), 398--404.

[Ja1] R.Janz, Holomorphic families of subspaces of a Banach space,
Special classes of linear operators and other  topics  (Bucharest,
1986), 155--167, Oper. Theory: Adv. Appl., 28,  Birkhauser,
Basel-Boston, MA, 1988.

[Ja2] R.Janz, Pertubation of Banach  spaces,  preprint,  Konstanz,
1987.

[JR] W.B.Johnson and  H.P.Rosenthal,  On $w^{*}$-basic  sequences  and
their applications to the study of Banach spaces, Studia  Math.,
43 (1972), 77--92.

[KT] W.Kaballo and G.Ph.A.Thijsse, On holomorphic operator  function
equations, Integral  Equations  and  Operator  Theory,  2  (1979),
no. 2, 244--263.

[Kad1] M.I.Kadets, On the linear dimension of the spaces
$L_p$  and
$l_q$ ,
Uspekhi Matem. Nauk, 13 (1958), no. 6, 95--98 (in Russian).

[Kad2] M.I.Kadets, Note on the gap between subspaces,  (in  Russian)
Funkts. Anal. Prilozhen. 9 (1975), no. 2, 73--74; English  transl.:
Funct. Anal. Appl. 9 (1975), 156--157.

[K1] T.Kato, Perturbation theory for nullity, deficiency and other
quantities of linear operators,  J. d'Analyse  Math.  6  (1958),
no. 2, 261--322.

[K2] T.Kato, Perturbation theory  for  linear  operators,  Berlin,
Springer-Verlag, 1966.

[KP] N.J.Kalton and N.T.Peck, Twisted sums of sequence spaces  and
the three space problem, Trans. Amer.  Math.  Soc.  255  (1979),
1--30.

[Ko1] C.A.Kottman,  Packing  and  reflexivity  in  Banach  spaces,
Trans. Amer. Math. Soc., 150 (1970), no. 2, 565--576.

[Ko2] C.A.Kottman, Subsets of the unit ball that are separated  by
more than one, Studia Math. 53 (1975), 15--27.

[KK] M.G.Krein and M.A.Krasnoselskii,  Basic  theorems  concerning
the  extension  of  Hermitian  operators  and  some   of   their
applications to the theory of  orthogonal  polynomials  and  the
moment problem, Uspekhi Mat. Nauk, 2 (1947), no. 3,  60--106  (in
Russian).

[KKM] M.G.Krein, M.A.Krasnoselskii and D.P.Milman, On  the  defect
numbers of linear operators  in  a  Banach  space  and  on  some
geometric questions, Sbornik Trudov Inst.  Matem.  AN  Ukrainian
SSR, 11 (1948), 97--112 (in Russian).

[KS] M.Krein and V.Smulian, On regularly closed sets in the  space
conjugate to a Banach space,  Ann.  Math.,  41  (1940),  no.  3,
556--583.

[Kw] S.Kwapien,  Isomorphic  characterizations  of  inner  product
spaces by orthogonal series  with  vector  valued  coefficients,
Studia Math., 44 (1972), no. 6, 583--595.

[L1] J.-P.Labrousse, On a metric space  of  closed  operators on a
Hilbert  space,  Rev.  Mat.  Fis.  Teorica  Univ.   N.   Tucuman
(Argentine) XVI (1 et 2) (1966), 45--77.

[L2]  J.-P.Labrousse,  Quelques   topologies   sur   des   espaces
d'operateurs dans des espaces de Hilbert  et  leurs  applications,
Dept. de Math. Univ. de Nice (1970).

[L3] J.-P.Labrousse, Metric on operator  spaces,  Lecture  given  at
Conference in Jaca (Spain), September, 1992.

[LM]  J.-P.Labrousse  and  B.Mercier,  Extension du   theoreme   de
Brown - Douglas -
Fillmore au  cas  des  operateurs  non  bornes,  J.
Operator Theory, 24 (1990), 137--154.

[Le] A.Yu.Levin, On a problem of existence of an orthogonal  element
to a subspace, Proc. Semin.  Funct.  Anal.  (Voronezh)  6  (1958),
91--92 (in Russian).

[LT1] J.Lindenstrauss and  L.Tzafriri,  Classical  Banach  spaces,
v. I, Sequence spaces, Berlin, Spriger-Verlag, 1977.

[LT2] J.Lindenstrauss and  J.Tzafriri,  Classical  Banach  spaces,
v. II, Function spaces, Berlin, Springer-Verlag, 1979.

[Lo1] W.E.Longstaff, A note on transforms of subspaces  of   Hilbert
space, Proc. Amer. Math. Soc., 76 (1979), no. 2, 268--270.

[Lo2] W.E.Longstaff, Subspace maps of operators  on  Hilbert  space,
Proc. Amer. Math. Soc., 84 (1982), no. 2, 195--201.

[Lu] W.Lusky, A note  on  Banach  spaces  containing $c_0$  and
$C_{\infty}$ , J.
Funct. Anal., 62 (1985), no. 1, 1--7.

[Man] F.Mantlik, Differentiable bundles of subspaces  and  operators
in Banach spaces, Studia Math. 98 (1991), 41--51.

[M1] A.S.Markus, On some properties of linear operators  connected
with the notion of opening, Uchen. Zap. Kishinev Gos.  Univ.  39
(1959), 265--272 (in Russian).

[M2] A.S.Markus, On comparability of openings in a  normed  space,
Izv. AN Moldov. SSR (1963), no. 1, 75--76 (in Russian).

[MS1] J.L.Massera and J.J.Sch\"affer, Linear differential  equations
and functional analysis, Ann. Math., 67 (1958), no. 3, 517--573.

[MS2] J.L.Massera and J.J.Sch\"affer, Linear Differential  Equations
and Function Spaces, New York and London, Academic Press, 1966.

[Ma] B.Maurey, Un theoreme de prolongement, C. r. Acad. Sc.  Paris,
Ser. A, 279 (1974), no. 9, 329--332.

[MeS] R.Mennicken and B.Sagraloff, On Banach's closed range theorem,
Arch. Math. 33 (1979/80), no. 5, 461--465.

[MM] D.P.Milman and V.D.Milman, The geometry of nested families with
empty intersection. Structure of the unit sphere of a nonreflexive
space, Matem. Sbornik, 66 (1965), no.  1,  109--118  (in  Russian).
English transl.: Amer.  Math.  Soc.  Transl.  (2)  v.~85  (1969),
233--243.

[MiS] V.D.Milman and G.Schechtman,  Asymptotic  Theory  of  Finite
Dimensional Normed Spaces, Berlin, Springer-Verlag, 1986.

[Na] M.A.Naimark, Linear Differential  Operators,  Moscow,  Nauka,
1969 (in Russian).

[N1] G.Neubauer,  \"Uber  den  Index  abgeschlossener Operatoren  in
Banachra\"umen, Math. Ann. 160 (1965), no. 2, 93--130.

[N2] G.Neubauer, \"Uber  den  Index  abgeschlossener  Operatoren  in
Banachra\"umen, II, Math. Ann. 162 (1965),no. 1, 92--119.

[New1] J.D.Newburgh, A topology for closed operators, Ann. Math., 53
(1951), 249--255.

[New2] J.D.Newburgh, The variation of  spectra,  Duke  Math.  J.  18
(1951), 165--176.

[O1] M.I.Ostrovskii, On the properties of the opening and  related
closeness characterizations of Banach spaces (in Russian), Teor.
Funktsii, Funktsional. Anal. i Prilozhen.,  42  (1984),  97--107:
Enslish translation:  Amer.  Math.  Soc.  Transl.  (2),  v.  136
(1987), 109--119.

[O2] M.I.Ostrovskii, The relation of  the  gaps  of  subspaces  to
dimension and other isomorphic invariants,  Akad.  Nauk  Armyan.
SSR Dokl. 79 (1984), no. 2, 51--53 (in Russian).

[O3] M.I.Ostrovskii, Banach-Saks properties, injectivity and  gaps
between subspaces of a Banach space (in Russian) Teor. Funktsii,
Funktsional. Anal.  i  Prilozhen.,  44  (1985),  69--78;  English
translation: J. Soviet Math., 48 (1990), no. 3, 299--306.

[O4] M.I.Ostrovskii, Closeness characterizations of Banach  spaces
generated by  the  notion  of  the  opening  between  subspaces,
Dissertation, Kharkov, 1985 (in Russian).

[O5] M.I.Ostrovskii, Separably injective Banach spaces (in Russian),
Funktsional. Anal. i Prilozhen. 20 (1986), no. 2, 80--81; English
translation.: Functional Anal. Appl. 20 (1986), no. 2, 154--155;

[O6] M.I.Ostrovskii, Properties of Banach spaces stable and unstable
with respect to the  gap  (in  Russian),  Sibirsk.  Mat.  Zh.  28
(1987), no. 1, 182--184; English translation: Siberian Math. J.,
28  (1987), no. 1, 140--142.

[O7] M.I.Ostrovskii, Comparison of topologies generated by geometric
and operator gaps between subspaces (in Russian), Ukr. Mat. Zh.
41 (1989), no. 7, 929--933; English translation: Ukr. Math. J.  41
(1989), no. 7, 794--797.

[O8] M.I.Ostrovskii, The three space problem for the weak Banach  -
Saks property  (in  Russian)  Mat.  Zametki  38  (1985),  no.  5,
721--725; English translation: Math. Notes  38  (1985),  no.  5--6,
905--907;

[O9] M.I.Ostrovskii, Classes of Banach spaces  stable  and  unstable
with respect to the opening, preprint.

[OP] M.I.Ostrovskii and A.N.Plichko, Banach-Saks properties  and the
three  space  problem,  in: "Operators  in  Function  Spaces  and
Problems in Function Theory", ed. by V.A.Marchenko, Kiev, Naukova
Dumka, 1987, p. 96--105 (in Russian).

[Par] K.R.Partasarathy, Lectures on functional  analysis.  Part  II.
Perturbations  by  unbounded  operators,  ISI  Lecture  Notes   6,
Macmillan Co of India Ltd, New Delhi, 1979.

[Pel1] A.Pelczynski, Projections in certain Banach spaces,  Studia
Math., 19 (1960), 209--228.

[Pel2] A.Pelczynski, Linear  extensions,  linear  averagings,  and
their  applications  to  linear  topological  classification  of
spaces of continuous functions, Warszawa, Panstwowe  Wydawnictwo
Naukowe, 1968.

[Pel3] A.Pelczynski, All separable Banach spaces admit for every
$\varepsilon >0$
fundamental total  and  bounded  by $1+\varepsilon $  biorthogonal  sequences,
Studia. Math. 55 (1976), 295--304.

[Pis1] G.Pisier, Factorization of Linear Operators and  Geometry  of
Banach Spaces, Providence, AMS, 1986.

[Pis2] G.Pisier, Quotients of Banach spaces of cotype $q$, Proc. Amer.
Math. Soc. 85 (1982), no. 1, 32--36.

[PX] G.Pisier and Q.Xu, Random series in real  interpolation  spaces
between the spaces $v_p$ , Lect. Notes Math. 1267 (1987), 185--209.

[Pl] A.Plichko, On bounded biorthogonal  systems  in  some  function
spaces, Studia Math. 84 (1986), 25--37.

[PR] H.Porta and L.Recht, Continuous selections of complemented sub%
spaces.  Contemporary Math. 52 (1986), 121--125

[Rak]  S.A.Rakov,  Ultraproducts  and  the  "three  space  problem",
Funktsional. Anal. i Prilozhen.,  11  (1977),  no.  3,  88--89  (in
Russian).

[RS] M.Reed and B.Simon, Methods of Modern Mathematical Physics,  1,
Functional Analysis, New York and London, Academic Press, 1972.

[Rod] L.Rodman, On global geometric properties  of  subspaces  in
Hilbert space, J. Funct. Anal. 45 (1982), 226--235.

[Ros] H.P.Rosenthal,  On  the  subspaces  of $L^{p} (p>2)$  spanned  by
sequences of independent random  variables,  Israel  J.  Math.,  8
(1970), no. 3, 273--303.

[Sc]  J.J.Sch\"affer,   Girth,   super-reflexivity,   and   isomorphic
classification of normed spaces and subspaces, Bull. Acad.  Polon.
Sci., Ser. Math., 28 (1980), no. 11--12, 573--584.

[ScS] J.J.Sch\"affer and K.Sundaresan, Reflexivity and  the  girth  of
spheres, Math. Ann., 184 (1970), 163--168.

[S1] I.Singer, Bases in Banach spaces, I,  Berlin,  Springer-Verlag,
1970.

[S2] I.Singer, Bases in Banach spaces, II, Berlin,  Springer-Verlag,
1981.

[S3]  I.Singer,  Best  approximation  in  normed  linear  spaces  by
elements of linear subspaces, Berlin, Springer-Verlag, 1970.

[Sl] Z.Slodkowski, Operators with closed ranges in  spaces  of
analytic vector-valued functions, J. Funct. Anal. 69 (1986), 155--177.

[Sob] P.E.Sobolevskii, On equations   with   operators   forming  an
acute angle, Dokl.  AN  SSSR,  116  (1957),  no.  5,  754--757  (in
Russian).

[Sz] S.J.Szarek, On the  geometry  of  the  Banach-Mazur  compactum,
Lect. Notes Math., 1470 (1991), 48--59.

[SzT] S.J.Szarek and M.Talagrand, An  "isomorphic"  version  of  the
Sauer-Shelah lemma and the  Banach-Mazur  distance  to  the  cube,
Lect. Notes Math., 1376 (1989), 105--112.

[Sz1] B. Sz.-Nagy, Perturbations  des  transformations autoadjointes
dans  l'espace  de  Hilbert,  Comment.  Math.  Helv.,  19  (1947),
347--366.

[Sz2]  B.  Sz.-Nagy,  Perturbations  des transformations   lineaires
fermees, Acta Sci. Math., 14 (1951), no. 2, 125--137.

[Te]  P.Terenzi,  On  the  theory  of  fundamental  norming  bounded
biorthogonal systems in Banach spaces, Trans. Amer. Math. Soc. 299
(1987), 497--511.

[T] V.M.Tikhomirov, Widths of sets in function spaces and the theory
of best approximations, Uspekhi  Mat.  Nauk,  15  (1960),  no.  3,
81--120, English translation: Russian  Math.  Surveys,  15  (1960),
no. 3, 75--111.

[V] M.Valdivia, Banach spaces $X$ with $X^{**}$ separable, Israel J.  Math.
59 (1987), no. 1, 107--111.

[Va1] F.-H.Vasilescu, Stability of the index of a complex of  Banach
spaces, J. Operator Theory, 2 (1979), no. 2, 247--275.

[Va2] F.-H.Vasilescu, Nonlinear  objects  in  the  linear  analysis,
Spectral Theory of Linear Operators and Related Topics (Timisoara/
Herculane, 1983), 265--278, Operator Theory, Adv. Appl., 14,  Birk%
hauser, Basel-Boston, Mass, 1984.

[Va3] F.~--H.~Vasilescu, Homogeneous operators and essential
complexes,
Glasgow
Math. J. 31 (1989), no. 1, 73--85.

[WW]  J.R.Wells  and  L.R.Williams,  Embeddings  and  Extensions  in
Analysis, Berlin, Springer-Verlag, 1975.

\end{large}
\end{document}